# Arithmetical Tugs of War and Benford's Law

**ABSTRACT**

Benford's Law predicts that the first significant digit on the leftmost side of numbers in real-life data is proportioned between all possible 1 to 9 digits approximately as in LOG(1 + 1/digit), so that low digits occur much more frequently than high digits in the first place. The two essential prerequisites for data configuration with regards to compliance with Benford's Law are high order of magnitude and positive skewness with a tail falling to the right of the histogram, so that quantitative configuration is such that the small is numerous and the big is rare. A related topic in the study of Benford's Law is the stark contrast between multiplications and additions of random variables and their distinct resultant quantitative and digital configurations. Random multiplication processes induce substantial increase in order of magnitude and they tend to the skewed Lognormal Distribution, favoring the small over the big. Random addition processes on the other hand do not induce any increase in order of magnitude and they tend to the symmetrical Normal Distribution as predicated by the Central Limit Theorem, favoring the medium over the small and the big. Thus, while multiplication processes are highly conducive to Benford behavior, addition processes are highly detrimental to Benford behavior. In this article it is shown that often in real-life data, multiplication and addition processes mix together within one measurement or expression, and consequently they fiercely compete for dominance, each attempting to exert the greatest influence upon sizes and digits. Such tugs of war between additions and multiplications are won or lost depending on the orders of magnitude of the generating random variables, as well as on the relative strength of the two warring sides, measured in terms of the comparative arithmetical involvement in the algebraic expression of the process.

Alex Ely Kossovsky
akossovsky@gmail.com

# PART 1: BENFORD'S LAW

# PART 2: ARITHMETICAL TUGS OF WAR

# PART 1:
# BENFORD'S LAW

# [I] The First Digit on the Left Side of Numbers

It has been discovered that the first digit on the left-most side of numbers in real-life data sets is most commonly of low value such as {1, 2, 3} and rarely of high value such as {7, 8, 9}. As an example serving as a brief and informal empirical test, a small sample of 40 values relating to geological data on time between earthquakes is randomly chosen from the data set on all global earthquake occurrences in 2012 – in units of seconds. Figure A depicts this small sample of 40 numbers. Figure B emphasizes in bold and black color the 1st digits of these 40 numbers.

| 285.29 | 185.35 | 2579.80 | 27.11 |
|--------|--------|---------|------|
| 5330.22 | 1504.49 | 1764.41 | 574.46 |
| 1722.16 | 815.06 | 3686.84 | 1501.61 |
| 494.17 | 362.48 | 1388.13 | 1817.27 |
| 3516.80 | 5049.66 | 2414.06 | 387.78 |
| 4385.23 | 2443.98 | 2204.12 | 1224.42 |
| 1965.46 | 3.61 | 1347.30 | 271.23 |
| 3247.99 | 753.80 | 1781.45 | 593.59 |
| 1482.64 | 1165.04 | 4647.39 | 1219.19 |
| 251.12 | 7345.52 | 1368.79 | 4112.13 |

**Figure A**: Sample of 40 Time Intervals between Earthquakes

| **2**85.29 | **1**85.35 | **2**579.80 | **2**7.11 |
|--------|--------|---------|------|
| **5**330.22 | **1**504.49 | **1**764.41 | **5**74.46 |
| **1**722.16 | **8**15.06 | **3**686.84 | **1**501.61 |
| **4**94.17 | **3**62.48 | **1**388.13 | **1**817.27 |
| **3**516.80 | **5**049.66 | **2**414.06 | **3**87.78 |
| **4**385.23 | **2**443.98 | **2**204.12 | **1**224.42 |
| **1**965.46 | **3**.61 | **1**347.30 | **2**71.23 |
| **3**247.99 | **7**53.80 | **1**781.45 | **5**93.59 |
| **1**482.64 | **1**165.04 | **4**647.39 | **1**219.19 |
| **2**51.12 | **7**345.52 | **1**368.79 | **4**112.13 |

**Figure B**: The First Digits of the Earthquake Sample

Clearly, for this very small sample, low digits occur by far more frequently on the first position than do high digits. A summary of the digital configuration of the sample is given as follows:

Digit Index:                                { 1,  2,  3,  4,  5, 6, 7, 8, 9 }
Digits Count totaling 40 values:            { 15, 8,  6,  4,  4, 0, 2, 1, 0 }
Proportions of Digits with '%' sign omitted: {38, 20, 15, 10, 10, 0, 5, 3, 0 }

Assuming (correctly) that these 40 values were collected in a truly random fashion from the large data set of all 19,452 earthquakes occurrences in 2012; without any bias or attempt to influence first digits occurrences; and that this pattern is generally found in many other data sets, one then may conclude with the phrase "not all digits are created equal", or rather "not all first digits are created equal", even though this seems to be contrary to intuition and against all common sense.

The focus here is actually on the first meaningful digit – counting from the left side of numbers, excluding any possible encounters of zero digits which only signify ignored exponents in the relevant set of powers of ten of our number system. Therefore, the complete definition of the **First Leading Digit** is the first non-zero digit of any given number on its left-most side. This digit is the first significant one in the number as focus moves from the left-most position towards the right, encountering the first non-zero digit signifying some quantity; hence it is also called the **First Significant Digit**. For 2365 the first leading digit is 2. For 0.00913 the first leading digit is 9 and the zeros are discarded; hence even though strictly-speaking the first digit on the left-most side of 0.00913 is 0, yet, the first significant digit is 9. For the lone integer 8 the leading digit is simply 8. For negative numbers the negative sign is discarded, hence for -715.9 the leading digit is 7. Here are some more illustrative examples:

**6**,719,525    → digit 6
0.0000**7**61    → digit 7
-0.**2**81264    → digit 2
**8**75          → digit 8
**3**            → digit 3
-**5**           → digit 5

For a data set where all the values are greater than or equal to 1, such as in the sample of the earthquaqe data, the first digit on the left-most side of numbers is also the First Leading Digit and the First Significant Digit, and necesarily one of the nine digits {1, 2, 3, 4, 5, 6, 7, 8, 9}; while digit 0 never occurs first on the left-most side.

# [II]   Benford's Law and the Predominance of Low Digits

Benford's Law states that:

Probability[First Leading Digit is d]  =  $LOG_{10}(1 + 1/d)$

$LOG_{10}(1 + 1/1) = LOG(2.00)$  = 0.301
$LOG_{10}(1 + 1/2) = LOG(1.50)$  = 0.176
$LOG_{10}(1 + 1/3) = LOG(1.33)$  = 0.125
$LOG_{10}(1 + 1/4) = LOG(1.25)$  = 0.097
$LOG_{10}(1 + 1/5) = LOG(1.20)$  = 0.079
$LOG_{10}(1 + 1/6) = LOG(1.17)$  = 0.067
$LOG_{10}(1 + 1/7) = LOG(1.14)$  = 0.058
$LOG_{10}(1 + 1/8) = LOG(1.13)$  = 0.051
$LOG_{10}(1 + 1/9) = LOG(1.11)$  = 0.046
                                  ---------
                                    1.000

Figure C depicts the distribution. Figure D visually depicts Benford's Law as a bar chart. This set of nine proportions of Benford's Law is sometimes referred to in the literature as **'The Logarithmic Distribution'**. Remarkably, Benford's Law is confirmed in almost all real-life data sets with high order of magnitude, such as in data relating to physics, chemistry, astronomy, economics, finance, accounting, geology, biology, engineering, governmental census data, and many others.

| Digit | Probability |
|---|---|
| 1 | 30.1% |
| 2 | 17.6% |
| 3 | 12.5% |
| 4 | 9.7% |
| 5 | 7.9% |
| 6 | 6.7% |
| 7 | 5.8% |
| 8 | 5.1% |
| 9 | 4.6% |

**Figure C**:  Benford's Law for First Digits

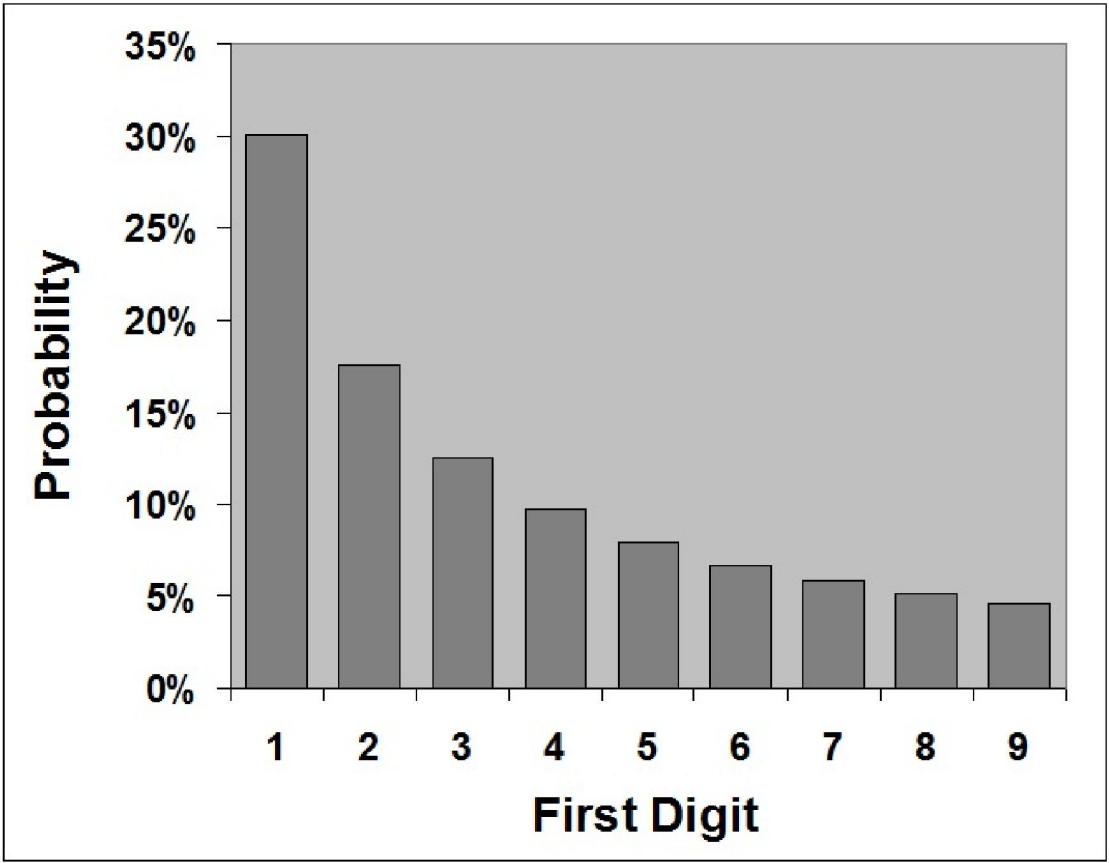

**Figure D**: Benford's Law – Probability of First Leading Digit Occurrences as a Bar Chart

# [III] Sum of Squares Deviation Measure (SSD)

It is necessary to establish a standard measure of 'distance' from the Benford digital configuration for any given data set. Such a numerical measure could perhaps tell us about the conformance or divergence from the Benford digital configuration of the data set under consideration. This is accomplished with what is called **Sum Squares Deviations (SSD)** defined as the sum of the squares of the 'errors' between the Benford expectations and the actual/observed values (in percent format – as opposed to fractional/proportional format):

$$SSD = \sum_{1}^{9} \left( \text{Observed \% of d} - 100 * \text{LOG}\left(1 + \frac{1}{d}\right) \right)^2$$

For example, for observed 1st digits proportions of {31.1, 18.2, 13.3, 9.4, 7.2, 6.3, 5.9, 4.5, 4.1} with '%' sign omitted, SSD measure of distance from the logarithmic is calculated as:

**SSD** = $(31.1 - \mathbf{30.1})^2 + (18.2 - \mathbf{17.6})^2 + (13.3 - \mathbf{12.5})^2 + (9.4 - \mathbf{9.7})^2 +$
$+ (7.2 - \mathbf{7.9})^2 + (6.3 - \mathbf{6.7})^2 + (5.9 - \mathbf{5.8})^2 + (4.5 - \mathbf{5.1})^2 + (4.1 - \mathbf{4.6})^2 = \mathbf{3.4}$

SSD generally should be below 25; a data set with SSD over 100 is considered to deviate too much from Benford; and a reading below 2 is considered to be ideally Benford.

# [IV]  Physical Order of Magnitude of Data (POM)

Rules regarding expectations of compliance with Benford's Law rely heavily on measures of order of magnitude (i.e. variability) of data.

Physical order of magnitude of a given data set is a measure that expresses the extent of its variability. It is defined as the ratio of the maximum value to the minimum value. The data set is assumed to contain only positive numbers greater than zero.

Physical Order of Magnitude (POM)  = Maximum∕Minimum

The classic definition of order of magnitude involves also the application of the logarithm to the ratio maximum/minimum, transforming it into a smaller and more manageable number.

Order of Magnitude (OOM) = $LOG_{10}$(Maximum∕Minimum)
Order of Magnitude (OOM) = $LOG_{10}$(Maximum)  -  $LOG_{10}$(Minimum)

Since such logarithmic transformation has a monotonic one-to-one relationship with max/min, it does not provide for any new insight or information, but could rather be looked upon in a sense as simply the use of an alternative scale, still measuring the same thing. For this reason the complexity of logarithm can be avoided altogether by referring only to the simple POM measure.

The more profound reason for using **POM** instead of **OOM** is its feature as a universal measure of variability, totally independent of societal number system in use, as well as being independent on the arbitrary choice of base 10, derived from the chanced or random occurrence of us having 10 fingers. This is the motivation behind the use of the term 'physical', expressing real and physical measure of variability, divorced from any numerical inventions, and especially so when data relates to the natural world such as in scientific figures and physical information.

OOM is perhaps more appropriate for a single isolated number, where it is re-defined as simply $LOG_{10}$(Number) without any reference to maximum, minimum, or any ratio. If we can assume that that number is an integral whole number without any [trailing] fractional part, then an alternative meaning of this OOM definition is simply expressing how many digits approximately are necessary to write the number. Surely in general, the bigger the [integral] number the more digits it takes to write it! For example, $LOG_{10}$(8,200,135) = 6.9, which is about 7, and that's exactly how many digits the number involves. As another example, $LOG_{10}$(10,000,000) = 7.0 which is exactly one digit less than the number of digits involved in writing the number, namely 8 digits.

In the extreme case where all the numbers in the data set are identical, having the value R say, variability is then nonexistent, and POM = maximum/minimum =  R/R = 1.

NOTE: LOG(X) or log(x) notation in this article would always refer to our decimal base 10 number system, hence the more detailed notation of $LOG_{10}$(X) is often (but not always) avoided.

# [V]  Two Essential Requirements for Benford Behavior

One of the two essential prerequisites or conditions for data configuration with regards to compliance with Benford's Law is that the value of the order of magnitude of the data set should be approximately over 3; in other words, that $LOG_{10}$(Maximum/Minimum) > 3, and that therefore (Maximum/Minimum) > $10^3$. This in turn implies that the threshold POM value (separating compliance from non-compliance) is about 1000, namely that POM > 1000 constitutes the condition for compliance.

Skewness of data where the histogram comes with a prominent tail falling to the right is the second essential criterion necessary for Benford behavior. Indeed, most real-life physical data sets are generally skewed in the aggregate, so that overall their histograms have tails falling on the right, and consequently the quantitative configuration is such that the small is numerous and the big is rare, while low first digits decisively outnumber high first digits.

The <u>asymmetrical</u>, Exponential, Lognormal, k/x [and many other distributions] are typical examples of such quantitatively skewed configuration, and therefore they are approximately, nearly, or exactly Benford - respectively. The <u>symmetrical</u> Uniform, Normal, Triangular, Circular-like, and other such distributions are inherently non-Benford, or rather anti-Benford, as they lack skewness and do not exhibit any bias or preference towards the small and the low.

Symmetrical distributions are always non-Benford, no matter what values are assigned to their parameters. By definition they lack that asymmetrical tail falling to the right, and such lack of skewness precludes Benford behavior regardless of the value of their order of magnitude.  Order of magnitude simply does not play any role whatsoever in Benford behavior for symmetrical distributions. For example, first digits of the Normal($10^{35}$, $10^8$) or the Uniform(1, $10^{27}$) are not Benford at all, and this is so in spite of their extremely large orders of magnitude. In summary: Benford behavior in extreme generality can be found with the confluence of sufficiently large order of magnitude together with skewness of data - having a histogram falling to the right. The combination of skewness and large order of magnitude is not a guarantee of Benford behavior, but it is a strong indication of likely Benford behavior under the right conditions. Moderate [overall] quantitative skewness with a tail falling too gently to the right implies that digits are not as skewed as in the Benford configuration. Extreme [overall] quantitative skewness with a tail falling sharply to the right implies that digits are severely skewed, even more so than they are in the Benford configuration. The only one exception to the generic rule above requiring skewness as well as high order of magnitude for logarithmic behavior is the perfectly Benford k/x distribution defined over adjacent integral powers of ten such as (1, 10) having the very low order of magnitude value of 1.

Bowley Skewness for example, defined as [(Q3 – Q2) – (Q2 – Q1)] **/** [Q3 – Q1] is an intuitive measure of skewness but its numerical value fluctuates greatly across data sets. Calculated Bowley Skewness values for numerous logarithmic data sets and distributions do not yield any consistent result, except that all values come out above 0.3 and below 1.0, and which is consistent with the fact that all logarithmic data sets are positively skewed in the aggregate. In sharp contrast, almost all non-logarithmic data sets come out with decisively lower Bowley Skewness values below 0.25 and above 0. In contrast to Bowley's unstable value for logarithmic data sets, Benford's Law is a very consistent and almost exact measure of skewness, with very little fluctuations across logarithmic data sets.

**Related Log Conjecture** in Benford's Law states that whenever the density curve of the logarithm of the data starts from the bottom on the log-axis itself, rises up to some plateau or maximum point, and then falls back all the way down to the log-axis itself, mimicking in a sense an upside-down U-like curve, and where total span on the log-axis from the minimum log to the maximum log is at least 3 units, then these two essential prerequisites above are guaranteed to be satisfied, and in addition guaranteeing that the data itself is nearly or almost perfectly Benford. See Kossovsky (2014) chapter 63 for a full discussion and concrete examples about the conjecture.

## [VI]   Data Skewness is More Prevalent than Benford's Law

All data sets obeying Benford's Law (i.e. logarithmic data) are structured in such a way that in the aggregate there are more small quantities than big quantities. In other words, that in the aggregate the histogram is falling to the right, except perhaps in the beginning on the very left for low values where it temporarily rises for a very small portion of overall data, as well as in few and minors reversals along the way where it rises briefly. This quantitative configuration is called 'positive skewness' in mathematical statistics.

This highly prevalent phenomenon of skewness has by far much wider scope as it is much more common in the physical world and in the realm of abstract mathematics than the more particular Benford quantitative configuration. This statement does not imply that Benford's Law is not prevalent in scientific, physical, and numerous other data types, on the contrary, it is highly prevalent. The statement only implies that in almost all the counter examples and exceptions to Benford's Law, the phenomenon of skewness is still found, albeit with different quantitative configurations than that of the Benford one (and typically with milder skewness, but at times even skewer). This renders Benford's Law a subset of the more universal skewness phenomenon.

The assertion is derived from concrete experience with numerous real-life numerical examples and from general research in Benford's Law. While this discussion may sound vague, in fact it is rather a very essential overview in the entire quantitative phenomenon of Benford's Law and of real-life data analysis in general. For those statisticians and data analysts who have worked on data sets and the Benford phenomenon for many years, including doing theoretical research, this generic statement seems natural, fundamental, and quite necessary.

# PART 2: ARITHMETICAL TUGS OF WAR

# [1]  Multiplication Processes Lead to Skewness and Higher POM

Almost all random and deterministic multiplication processes induce a dramatic increase in skewness where the small becomes relatively numerous and the big becomes relatively rare. Surely, the single product of the multiplication of only two numbers is not to be considered here, since the resultant single product is neither small nor big, but rather it's just itself, and there are no competing sizes here to consider. Instead, a set of N numbers called A is to be multiplied by another set of M numbers called B. The phrase 'multiplication of two sets of numbers' implies that each and every number in set A is to be multiplied by each and every number of set B, producing N*M products. In other words, all possible multiplications between the two sets are attempted, and the entirety of these products constitutes the newly created set of numbers, called A*B. For example:

Set A = {8, 3, 5}
Set B = {11, 47, 26}
A*B = {(8)*(11), (8)*(47), (8)*(26), (3)*(11), (3)*(47), (3)*(26), (5)*(11), (5)*(47), (5)*(26)}
A*B = {88, 376, 208, 33, 141, 78, 55, 235, 130}

Let us examine the quantities within the 10 by 10 multiplication table that we all were forced to memorize in our elementary school years. In this example, the intrinsic characteristics of multiplication processes with regards to resultant sizes shall be explored and compared. Such an analysis is done at the most primitive and basic level, at the arithmetic and quantitative level, before going on to the rigorous mathematical level. The quest is to start out with this very particular example, and then to lend the conclusions derived from this case universality and applicability in almost all multiplication processes.

For the 10 by 10 multiplication table, the entire range from 1 to 100 is to be partitioned into 10 quantitative sections of equal 10-unit width each, namely [1, 10], [11, 20], [21, 30], [31, 40], [41, 50], [51, 60], [61, 70], [71, 80], [81, 90], [91, 100], and then a count is made of the values falling within each section. The goal is to group the 100 products of the multiplication table according to sizes. Figure 1 depicts this quantitative partitioning arrangement of the entire multiplicative territory by size. Figure 2 depicts the histogram of the values falling within each section. Surprisingly, a decisive trend regarding the occurrences of products within the sections is found here. The section of the smallest quantities (1 to 10) has 27 values falling within it; while the section of the biggest quantities (91 to 100) has only 1 value falling within it.

The sequence of all the values falling within these 10 sections, and presented in order according to sizes is {27, 19, 15, 11, 9, 6, 5, 4, 3, 1}. Clearly, the sections pertaining to bigger quantities have less values falling within them, as the count of values monotonically decreases. The small is definitely more numerous than the big in our standard multiplication table. The crucial lesson learnt from this multiplication process is that surely this tendency is nearly universal, and that it should be present in almost all other multiplication processes, and not only for our particular 10 by 10 multiplication table. There is nothing unique about our standard multiplication table or the set of numbers from 1 to 10, and so this result is extrapolated to almost all other multiplication processes.

The results of Figures 1 and 2 can also be interpreted as a particular casino game where two dice having 10 faces each are thrown, and the values of the two faces are multiplied by each other. The value of this product is then declared to be the focus of the gambling game, and an arbitrary benchmark value of say 30 is set determining losses and wins. Players are said to win if the product comes out in the range of 31 – 100 of bigger values, and the casino is said to win if the product comes out in the range of 1 – 30 of smaller values. Such cleaver setup on the part of the casino owner where the small is assigned to the casino and the big is assigned to the gamblers yields steady and enormous profits for the casino. Our 10 by 10 multiplication table can also be interpreted in a more formal sense in the context of mathematical statistics as the random multiplicative process of two **discrete** random Uniform Distributions, namely the product: Uniform{1, 2, 3, 4, 5, 6, 7, 8, 9, 10}*Uniform{1, 2, 3, 4, 5, 6, 7, 8, 9, 10}, instead of regarding it just as a useful tool of the deterministic table of multiplication. More generally, it could also be thought of as the random multiplicative process of two **continuous** random Uniform Distributions, namely the product Uniform[1, 11)*Uniform[1, 11).

| X | 1 | 2 | 3 | 4 | 5 | 6 | 7 | 8 | 9 | 10 |
|---|---|---|---|---|---|---|---|---|---|---|
| 1 | 1 | 2 | 3 | 4 | 5 | 6 | 7 | 8 | 9 | 10 |
| 2 | 2 | 4 | 6 | 8 | 10 | 12 | 14 | 16 | 18 | 20 |
| 3 | 3 | 6 | 9 | 12 | 15 | 18 | 21 | 24 | 27 | 30 |
| 4 | 4 | 8 | 12 | 16 | 20 | 24 | 28 | 32 | 36 | 40 |
| 5 | 5 | 10 | 15 | 20 | 25 | 30 | 35 | 40 | 45 | 50 |
| 6 | 6 | 12 | 18 | 24 | 30 | 36 | 42 | 48 | 54 | 60 |
| 7 | 7 | 14 | 21 | 28 | 35 | 42 | 49 | 56 | 63 | 70 |
| 8 | 8 | 16 | 24 | 32 | 40 | 48 | 56 | 64 | 72 | 80 |
| 9 | 9 | 18 | 27 | 36 | 45 | 54 | 63 | 72 | 81 | 90 |
| 10 | 10 | 20 | 30 | 40 | 50 | 60 | 70 | 80 | 90 | 100 |

**Figure 1**: Quantitative Territorial Partitioning of the 10 by 10 Multiplication Table

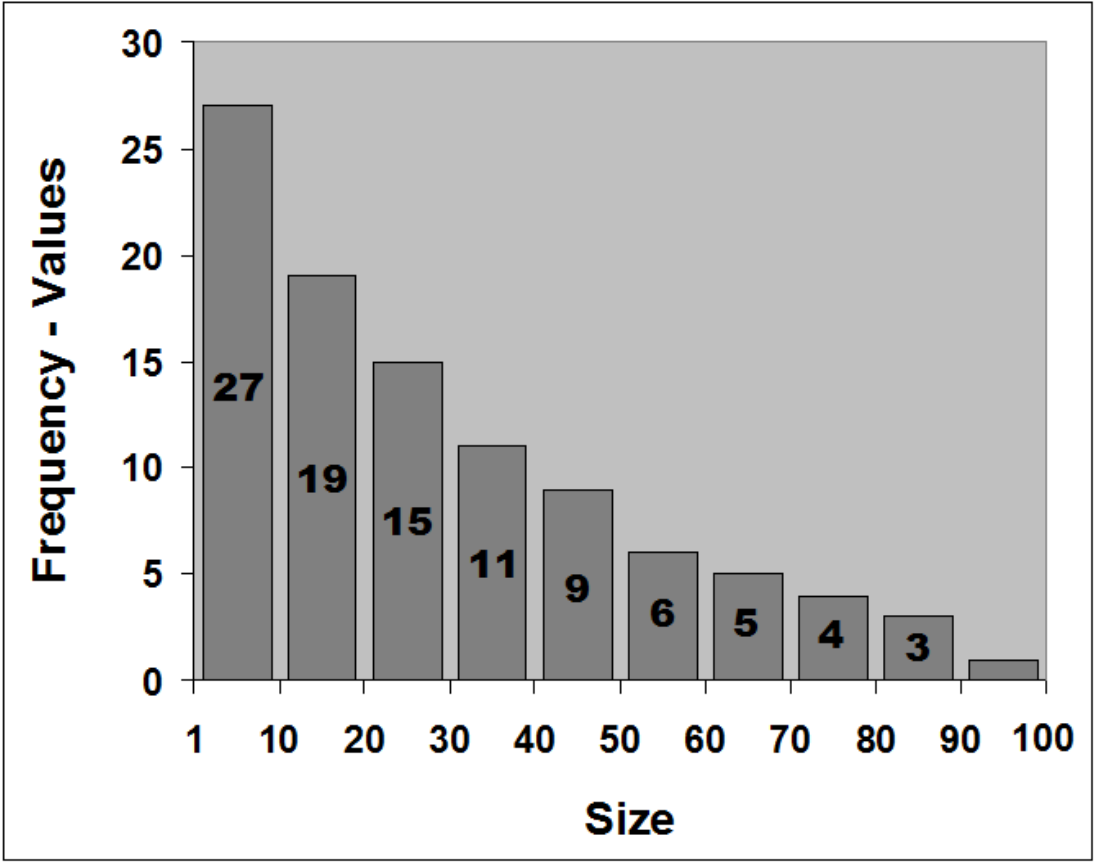

**Figure 2**: Histogram of the 10 by 10 Multiplication Table Depicting Acute Skewness

It should be noted that for the above result about the 10 by 10 multiplication table, the starting point is the uniform and even distribution of the original numbers about to be multiplied, namely {1, 2, 3, 4, 5, 6, 7, 8, 9, 10}, where neither small nor big is more numerous than the other. Yet, from such even and uniform distribution we arrived (via multiplications) at a decisively skewed and uneven distribution where the small is more numerous than the big.

Exception to this rule is found only in multiplications of rare sets of numbers with almost no variability, namely data sets having extremely low order of magnitude. This case shall be discussed in chapter 5 of this article.

In mathematical statistics, the random product of numerous independent realizations from an identical random variable is known to be the Lognormal Distribution, in the limit, as the number of these multiplied realizations becomes large enough. This seminal result which is derived from the Multiplicative Central Limit Theorem (MCLT) has very lax requirements and only a few restrictions, ensuring broad applicability for almost all types of multiplied variables, namely for almost all multiplication processes. Most significantly, the restriction on having an identical distribution can often be relaxed assuming other easily obtained conditions, thus MCLT can usually be applied to the random product of several distinct random variables. Since the Lognormal Distribution is almost always skewed [that is, for all shape parameters above 0.5 roughly], MCLT guarantees that almost all multiplication processes are skewed.

The scope of applications here is truly enormous, encompassing almost all disciplines. The customary multiplicative form of the vast majority of the equations in physics, chemistry, astronomy, economics, biology, engineering and other disciplines, as well as the numerous expressions relating to their applications and specific results, leads to the manifestation of quantitative skewness and often also to the Benford phenomenon in the natural sciences and in real-life data. These data sets are almost always with high enough variability (i.e. high order of magnitude) and thus they almost never suffer from the above-mentioned exception to the rule.

As discussed in Kossovsky (2014) chapter 90, Isaac Newton gave us $F = M*A$, not $F = M + A$. Hence the derivations and results due to the expression Force = Mass*Acceleration are related to multiplication processes, not to addition processes. Newton also gave us the law of universal gravitation $F_G = G*M_1*M_2 / R^2$ which is written in multiplicative and divisional forms, and not in addition and subtraction style such as say $F_G = G + M_1 + M_2 - R^2$.

Resultant quantitative skewness in multiplication processes often leads to digital skewness and Benford as well, hence strengthening the plausible explanation of Benford's Law in terms of multiplication processes and their intrinsic logarithmic properties.

For our 10 by 10 multiplication table from our childhood school days, order of magnitude doubles from 1 to 2 as we multiply; namely, from the range of [1, 10] we arrive at the range of [1, 100], so that variability dramatically increases due to the act of multiplying.
POM of [1, 10] is 10/1, and POM of [1, 100] is 100/1, thus there is a ten-fold increase in POM.
OOM of [1, 10] is LOG(10/1) = 1, and OOM of [1, 100] is LOG(100/1) = 2, thus OOM has been doubled due to the act of multiplying.

In summarizing the generic effects of random multiplication processes on quantities, it can be said that **multiplication processes favor the small over the big - leading to skewed data**.
In addition, when the focus is on digital configurations in the first position, it can be said that multiplication processes favor small digits such as {1, 2, 3} over big digits such as {7, 8, 9}.

Random multiplication processes induce two essential results:

**(A)** A dramatic increase in skewness – an essential criterion for Benford behavior.
**(B)** An increase in the order of magnitude – another essential criterion for Benford behavior.

# [2]  Addition Processes Lead to the Symmetrical Normal and Equal POM

Let us demonstrate the sharp contrast between additions and multiplications in terms of resultant quantitative configuration and resultant order of magnitude. This shall be accomplished by converting the 10 by 10 multiplication table into an addition table.

We may view the table as a tool for those too lazy to do additions quickly in their heads, just as the multiplication table is used and memorized as a tool to be used later in life. Or, we may view each discrete set of {1, 2, 3, 4, 5, 6, 7, 8, 9, 10} as a random variable with uniform and equal probability for each integer, and the table as an addition process. Since additions are applied instead of multiplications, such a random vista about this (essentially deterministic) addition table allows us perhaps to extrapolate the results to all random addition processes.

The entire range of [2, 20] is partitioned into six equitable quantitative sections of 3-unit width each: {2, 3, 4}, {5, 6, 7}, {8, 9, 10}, {11, 12, 13}, {14, 15, 16}, {17, 18, 19} *[with the value of 20 conveniently excluded]*. Then, a count is made of the numbers falling within each section - namely grouping them according to quantities. Figure 3 demonstrates such quantitative partitioning of the entire additive territory in details. Figure 4 depicts the histogram of the number of values falling within each section.

| + | 1 | 2 | 3 | 4 | 5 | 6 | 7 | 8 | 9 | 10 |
|---|---|---|---|---|---|---|---|---|---|---|
| 1 | 2 | 3 | 4 | 5 | 6 | 7 | 8 | 9 | 10 | 11 |
| 2 | 3 | 4 | 5 | 6 | 7 | 8 | 9 | 10 | 11 | 12 |
| 3 | 4 | 5 | 6 | 7 | 8 | 9 | 10 | 11 | 12 | 13 |
| 4 | 5 | 6 | 7 | 8 | 9 | 10 | 11 | 12 | 13 | 14 |
| 5 | 6 | 7 | 8 | 9 | 10 | 11 | 12 | 13 | 14 | 15 |
| 6 | 7 | 8 | 9 | 10 | 11 | 12 | 13 | 14 | 15 | 16 |
| 7 | 8 | 9 | 10 | 11 | 12 | 13 | 14 | 15 | 16 | 17 |
| 8 | 9 | 10 | 11 | 12 | 13 | 14 | 15 | 16 | 17 | 18 |
| 9 | 10 | 11 | 12 | 13 | 14 | 15 | 16 | 17 | 18 | 19 |
| 10 | 11 | 12 | 13 | 14 | 15 | 16 | 17 | 18 | 19 | 20 |

**Figure 3:**  Quantitative Territorial Partitioning of the Ten by Ten Addition Table

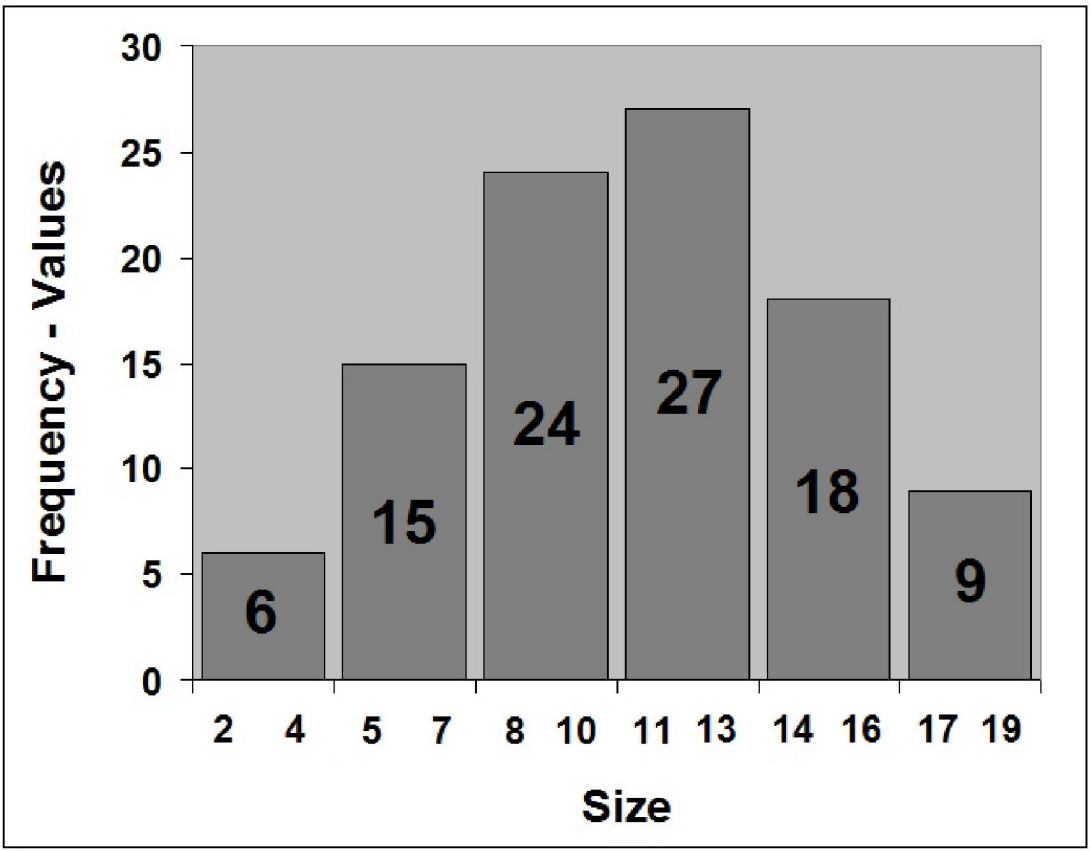

**Figure 4:** Histogram of Ten by Ten Addition Table Depicting Normal-Like Symmetry

There are  **6** values from  2 to  4.
There are **15** values from  5 to  7.
There are **24** values from  8 to 10.
There are **27** values from 11 to 13.
There are **18** values from 14 to 16.
There are  **9** values from 17 to 19.

The histogram is not skewed, and quantitative proportions are concentrated mostly around the middle part of the range, not around the left part (of small values), and not around the right part (of big values). This is not a coincidence, but rather it signifies a very persistent trend found in all addition processes. There is no hope that repeated additions of random variables would eventually attain skewness assuming large number of applied additions, because any such aspiration on the part of addition is frustrated by the Central Limit Theorem (CLT) which points to the symmetric Normal Distribution as the eventual distribution developing after numerous such additions of random variables. A strong hint of this eventuality can actually be seen in Figure 4 as it nicely curves around the center and falls off almost evenly on both edges, beginning to resemble the Normal Distribution.

Here for this 10 by 10 addition table, order of magnitude stubbornly refuses to increase. The maximum value is 20, and the minimum value is 2, hence resultant post-addition POM value for the entire table is (20)/(2) = 10. This value is exactly the same as the (10)/(1) = 10 POM value of the original variable being added, namely of the discrete set of {1, 2, 3, 4, 5, 6, 7, 8, 9, 10}.

Clearly, when the focus is on quantitative configurations, it can be said in extreme generality that **addition processes favor the medium over the small and over the big**. When the focus is on digital configurations in the first position, nothing in general can be said a priori about addition processes, and digital configuration depends on the specifics of the added variables and especially on defined ranges of the added variables. Addition processes do not favor a priori middle digits such as {4, 5, 6} except in some very particular cases.

Random addition processes do not induce any results that are essential in the criterion for Benford behavior:

**(A)** Lacking any increase in skewness, and even actively increasing the symmetry of resultant distribution, with added concentration forming around the center/medium.

**(B)** Lacking any increase in order of magnitude beyond the existing maximum order of magnitude within the set of added variables.

Addition processes are associated with mild and slow quantitative increases. For example, if we walk along the diagonal of the 10 by 10 addition table of Figure 3 *[from top-left of 1+1 to bottom-right of 10+10]*, the sequence 1+1, 2+2, 3+3, 4+4, 5+5, 6+6, 7+7, 8+8, 9+9, 10+10, increases steadily by 2 units at a time without any acceleration. Figure 5 depicts the steady and even march forward of the added ten numbers 2, 4, 6, 8, 10, 12, 14, 16, 18, 20, which are evenly spread along the horizontal axis.

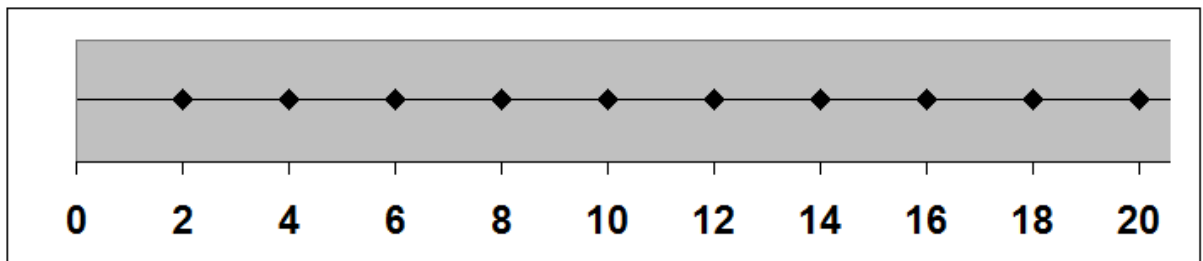

**Figure 5**: Even and Balanced Spread of the Added Numbers **X + X**, for X ϵ {1, 2, 3, … 10}

In sharp contrast, multiplication processes are associated with strong and rapid quantitative increases. Numbers under 'multiplicative influences' almost always undergo what could be called '**quantitative acceleration**'.

For example, if we walk along the diagonal of the squares in the 10 by 10 multiplication table of Figure 1 *[from top-left of 1*1 to bottom-right of 10*10]*, we experience a rapid rise of quantities, from 1*1, to 2*2, to 3*3, to 4*4, to 5*5, to 6*6, to 7*7, to 8*8, to 9*9, and to 10*10, pointing to the sequence of products of 1, 4, 9, 16, 25, 36, 49, 64, 81, 100. Figure 6 depicts the accelerated march forward of the square numbers 1 to 100. Walking along the diagonal of the squares from say **5*5** to **6*6**, means that now we are adding **6+6+6+6+6+6** instead of **5+5+5+5+5**, and which represents a quantitative jump in a sense, because not only an extra term is added, but also because each term is higher by one notch. The bigger the number, the more effective multiplication becomes in moving the product forward. In other words, the bigger the multiplicands, the (much) greater are the products.

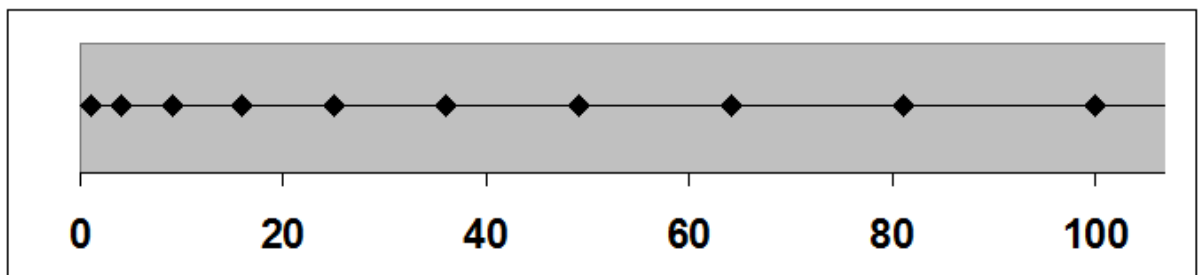

**Figure 6**: The Accelerated March of the Square Numbers **X*X**, for X ϵ {1, 2, 3, … 10}

# [3]  Arithmetically-Mixed Real-Life Variables

The 10 by 10 multiplication table involves purely multiplicative processes. The 10 by 10 addition table of the previous chapter involves a purely additive process. These are examples of exclusively multiplication processes or that of exclusively addition processes. Often in real-life data however, multiplication and addition processes mix together within one measurement, and consequently they fiercely compete for dominance, each attempting to exert the greatest influence upon sizes and digits. While addition favors the medium, multiplication prefers the small, and none is willing to compromise.

As an example relevant to real-life accounting and financial data, the single bill for one typical shopper in a large supermarket or at a big retail store may read as follows:

**3***($**2.75** bread) + **5***($**2.50** tuna) + **2***($**7.99** cheese) = ($**36.73** total bill)

Namely, buying 3 loafs of hearty brown whole wheat Russian bread; 5 cans of white Albacore tuna in olive oil; and 2 packages of low-fat and preservative-free Mozzarella cheese. Here multiplications mix in with additions, and as a consequence these two competing forces pull in opposite directions. One multiplicand within the expression of a typical bill is the catalog of a retail store or the price list of all the items on sale at a large supermarket. The other multiplicand is the quantity purchased of each item (such as the numbers 3, 5, 2 in the above bill of $36.73). The bills of thousands such potential buyers at the supermarket constitute a random variable derived from a mixture of addition and multiplication processes, representing revenue data for the supermarket and expense or cost data for the shopper. Surely, the expectation here is that the quantitative configuration of the catalog or the price list itself could greatly influence relative quantities and sizes of the bills.

As another example relevant to real-life scientific and chemical data, the weight of a complex chemical molecule is derived from the linear combination of its constituent atoms.

For example, Lactose $C_{12}H_{22}O_{11}$ has the molar mass of 342**.**29648 g/mol. This particular molecular weight is derived from the combinations:

**12***(Carbon Mass) + **22***(Hydrogen Mass) + **11***(Oxygen Mass) = Lactose Mass

Here one multiplicand is the atomic weight in the Periodic Table, and the other multiplicand is the number of atoms for each element within the molecule.

In all of these examples, additions and multiplications mix and combine to yield the final value, and each exerts its distinct pull and influence. What should be expected in much arithmetical mixtures in terms of (I) quantitative skewness versus symmetry, and in terms of (II) digital configuration?  Who wins and who loses in such tugs of war? Does it simply depend on the relative number of additions versus the relative number of multiplications within the arithmetical expression? The answers to these questions shall become clear in the next several chapters.

# [4] Arithmetically-Complex Random Processes

Let us examine the quantitative configuration of 4 arithmetically-complex dice games.

We view each standard 6-sided die as a discrete set of {1, 2, 3, 4, 5, 6}; consider it as a random variable with uniform and equal probability for each integer; and call it SIX.

Each of the four dice games shown below involves using 10 regular 6-sided dice combinations, as follows:

**Game 1**:  $(SIX_1*SIX_2*SIX_3) + (SIX_4*SIX_5*SIX_6) + (SIX_7*SIX_8*SIX_9) + (SIX_{10})$

**Game 2**:  $(SIX_1*SIX_2*SIX_3) + (SIX_4*SIX_5*SIX_6) + (SIX_7*SIX_8) + (SIX_9*SIX_{10})$

**Game 3**:  $(SIX_1*SIX_2*SIX_3) + SIX_4 + SIX_5 + SIX_6 + SIX_7 + SIX_8 + SIX_9 + SIX_{10}$

**Game 4**:  $(SIX_1*SIX_2) + (SIX_3*SIX_4) + (SIX_5*SIX_6) + (SIX_7*SIX_8) + (SIX_9*SIX_{10})$

Here each die has an identity and a unique name. They are designated $SIX_1$, $SIX_2$, $SIX_3$, and so forth up to $SIX_{10}$. The ten dice are thrown simultaneously and independently of each other, and then the arithmetically-complex operations are calculated on the 10 values appearing on the 10 faces of the 10 dice.

Figures 7, 8, 9, and 10 depict the histograms for these four games, using 10,000 Monte Carlo computer simulation runs per game.

For games 1, 2, and 3, of Figures 7, 8, and 9, the histogram appears as a hybrid of sorts, standing somewhere between the symmetric Normal histogram of additions, and the generally skewed histogram of multiplications. Perhaps they appear more skewed than symmetrical. This is expected, since these processes are truly mixtures of addition and multiplication operations.

For game 4 of Figure 10, the histogram is quite symmetrical, almost mimicking the Normal Distribution. Perhaps this is expected, and almost predicated by the Central Limit Theorem, since the process might be thought of primarily as that of additions. The term (SIX*SIX) is being added over and over again five times. Certainly the term (SIX*SIX) can be thought of as a singular random distribution in its own right, and thus being the addend in this process. The consideration of (SIX*SIX) as a whole entity leads to the simpler vista of having five random variables being added, and this eliminates the appearance of multiplication altogether from the process. In other words, such a vista facilitates the consideration of the final value in this game as derived from a purely addition process. A similar vista though could also be presented for games 1, 2, and 3, but since in those games distinct variables are being added, the Central Limit Theorem is slow to act, and therefore results are not as Normal-like.

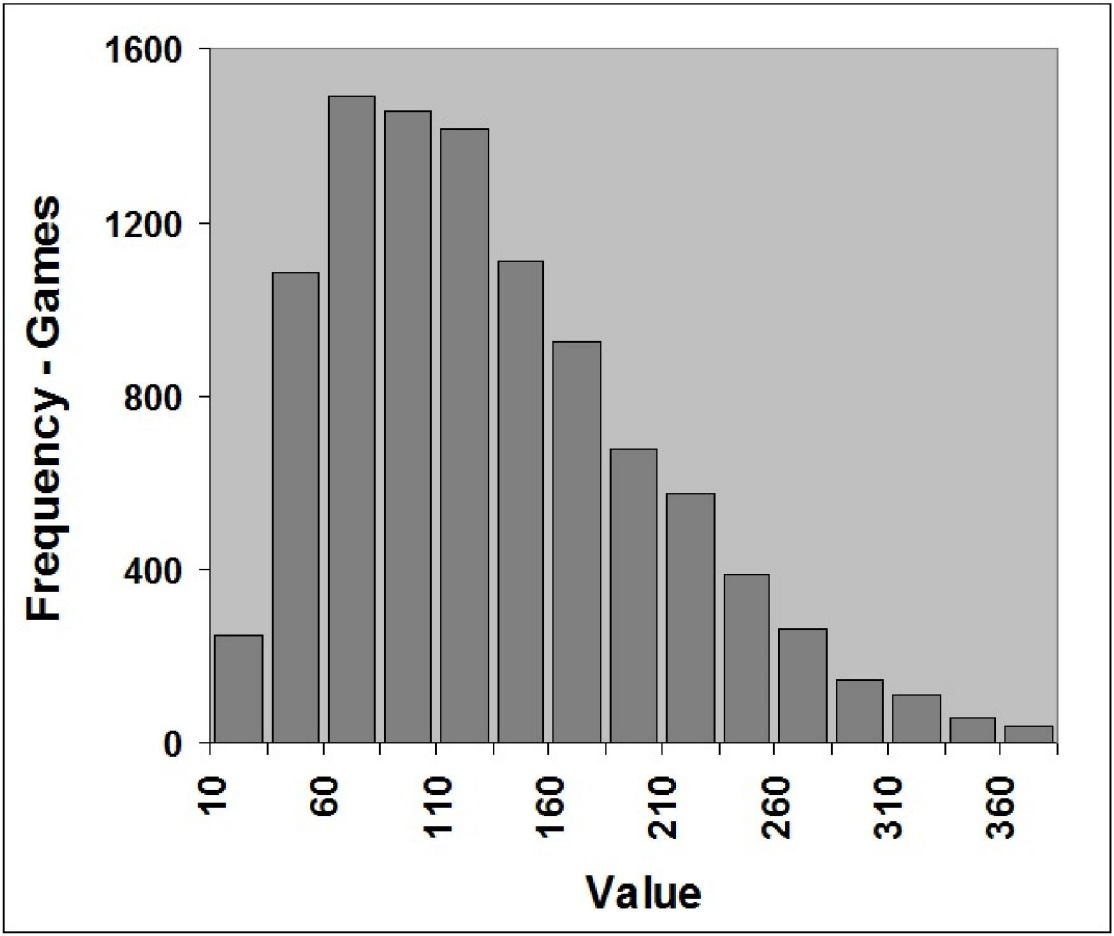

**Figure 7:** $(SIX_1*SIX_2*SIX_3) + (SIX_4*SIX_5*SIX_6) + (SIX_7*SIX_8*SIX_9) + (SIX_{10})$

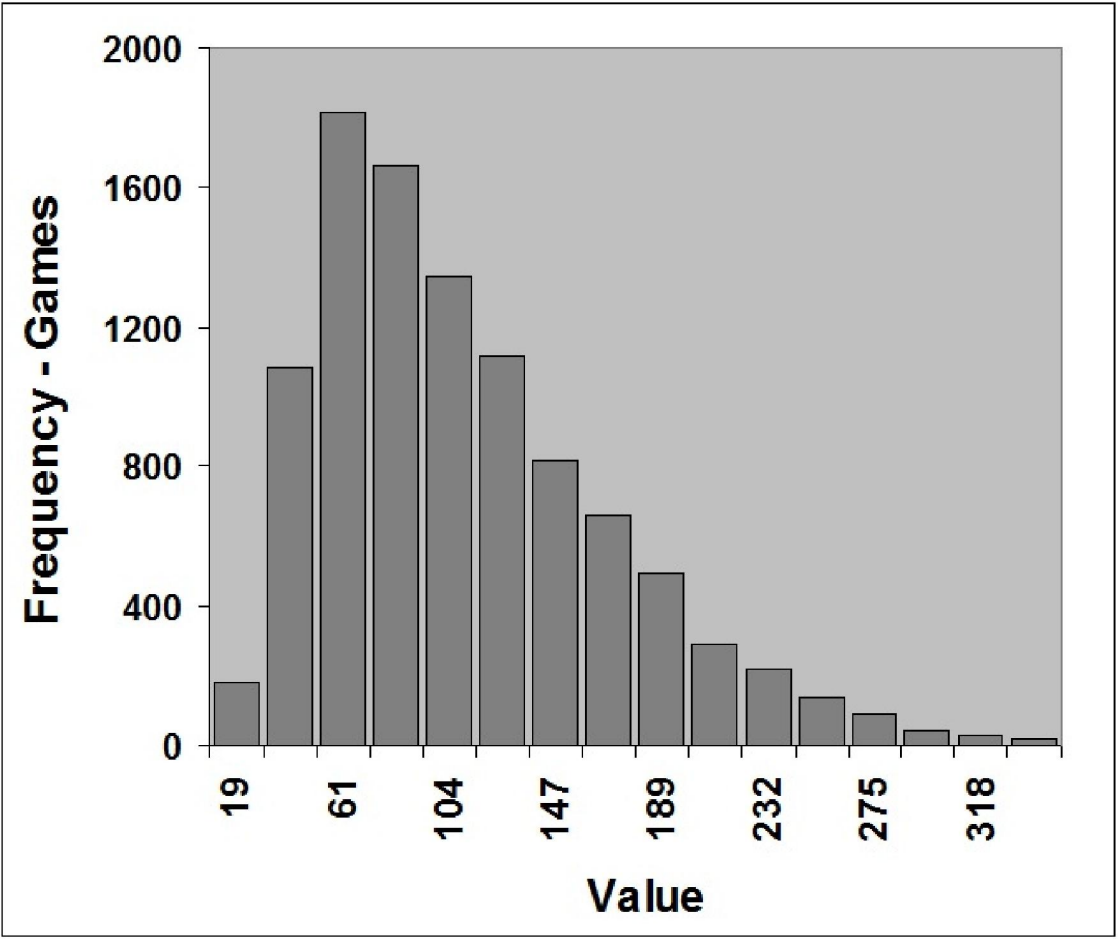

**Figure 8:** $(SIX_1*SIX_2*SIX_3) + (SIX_4*SIX_5*SIX_6) + (SIX_7*SIX_8) + (SIX_9*SIX_{10})$

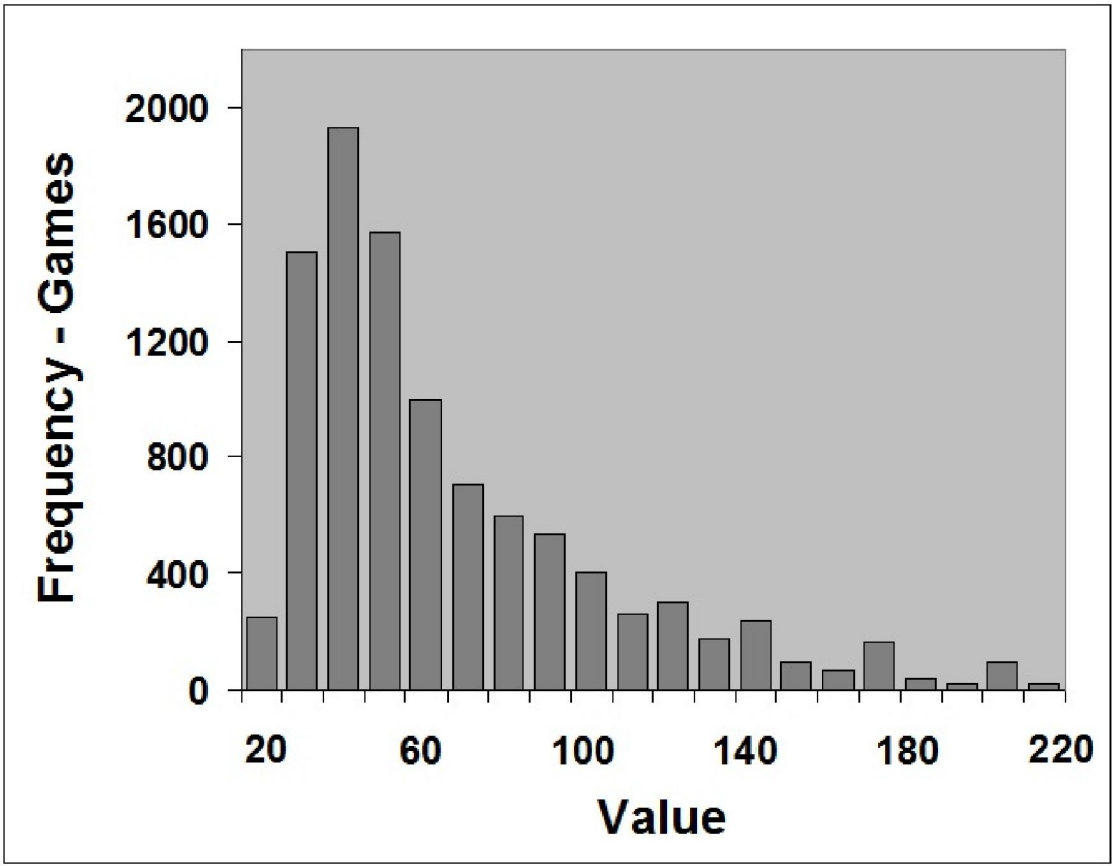

**Figure 9:** $(SIX_1 * SIX_2 * SIX_3) + SIX_4 + SIX_5 + SIX_6 + SIX_7 + SIX_8 + SIX_9 + SIX_{10}$

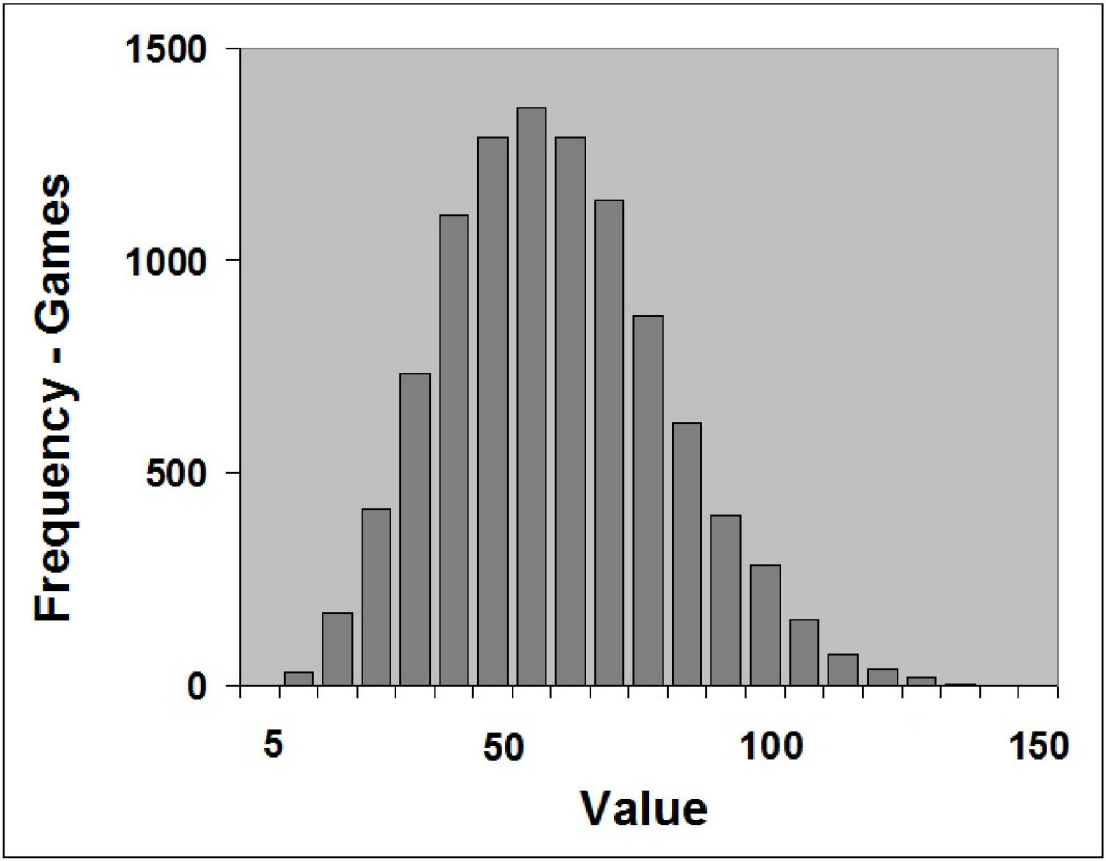

**Figure 10:** $(SIX_1*SIX_2) + (SIX_3*SIX_4) + (SIX_5*SIX_6) + (SIX_7*SIX_8) + (SIX_9*SIX_{10})$

### [5] Multiplications of Variables with Low Order of Magnitude

In order to enable a thorough analysis of the entire discourse regarding mixed arithmetical processes and their quantitative and digital properties, it is necessary to digress momentarily and explore random multiplication processes involving variables with low order of magnitude.

Let us examine quantities within the 70 to 79 multiplication table. This table is a variation on the standard 10 by 10 table, using the set of integers {70, 71, 72, 73, 74, 75, 76, 77, 78, 79} instead of {1, 2, 3, 4, 5, 6, 7, 8, 9, 10}. Figure 11 depicts the arrangement.

The motivation for choosing this set of integers is to examine skewness for multiplications of variables with very low physical order of magnitude. For this set of integers, POM (70 to 79) = Maximum/Minimum = (79)/(70) = 1.13. In contrast, for the standard 10 by 10 multiplication table, POM (1 to 10) = Maximum/Minimum = (10)/(1) = 10.0. Clearly, POM value for the 70 to 79 set of integers is extremely low, hovering just slightly above the lowest possible 1.0 value of POM for data without any variability whatsoever.

| x  | 70   | 71   | 72   | 73   | 74   | 75   | 76   | 77   | 78   | 79   |
|----|------|------|------|------|------|------|------|------|------|------|
| 70 | 4900 | 4970 | 5040 | 5110 | 5180 | 5250 | 5320 | 5390 | 5460 | 5530 |
| 71 | 4970 | 5041 | 5112 | 5183 | 5254 | 5325 | 5396 | 5467 | 5538 | 5609 |
| 72 | 5040 | 5112 | 5184 | 5256 | 5328 | 5400 | 5472 | 5544 | 5616 | 5688 |
| 73 | 5110 | 5183 | 5256 | 5329 | 5402 | 5475 | 5548 | 5621 | 5694 | 5767 |
| 74 | 5180 | 5254 | 5328 | 5402 | 5476 | 5550 | 5624 | 5698 | 5772 | 5846 |
| 75 | 5250 | 5325 | 5400 | 5475 | 5550 | 5625 | 5700 | 5775 | 5850 | 5925 |
| 76 | 5320 | 5396 | 5472 | 5548 | 5624 | 5700 | 5776 | 5852 | 5928 | 6004 |
| 77 | 5390 | 5467 | 5544 | 5621 | 5698 | 5775 | 5852 | 5929 | 6006 | 6083 |
| 78 | 5460 | 5538 | 5616 | 5694 | 5772 | 5850 | 5928 | 6006 | 6084 | 6162 |
| 79 | 5530 | 5609 | 5688 | 5767 | 5846 | 5925 | 6004 | 6083 | 6162 | 6241 |

**Figure 11**: Quantitative Territorial Partitioning of the 70 to 79 Multiplication Table

The lowest product in the 70 to 79 table is 70*70 = 4900, and the highest is 79*79 = 6241.
The intention is to partition the entire range from **4900** to **6241** into 10 quantitative sections of equal width and compare density within each section (i.e. count of data points falling in).
In order to obtain a smooth and integral partition not involving fractional values as border points, this range shall be widen a tiny bit to span **4900** to **6300**, and therefore the enlarged range is of (6300 – 4900) = 1400 units. Dividing 1400 into 10 equitable sections implies having the width of 140 units for each bin in the histogram.

Figure 12 depicts the histogram of the values falling within each section.

There are  3  values from  4900  up to  5040.
There are  7  values from  5040  up to  5180.
There are 11  values from  5180  up to  5320.
There are 15  values from  5320  up to  5460.
There are 19  values from  5460  up to  5600.
There are 17  values from  5600  up to  5740.
There are 13  values from  5740  up to  5880.
There are  9  values from  5880  up to  6020.
There are  3  values from  6020  up to  6160.
There are  3  values from  6160  up to  6300.

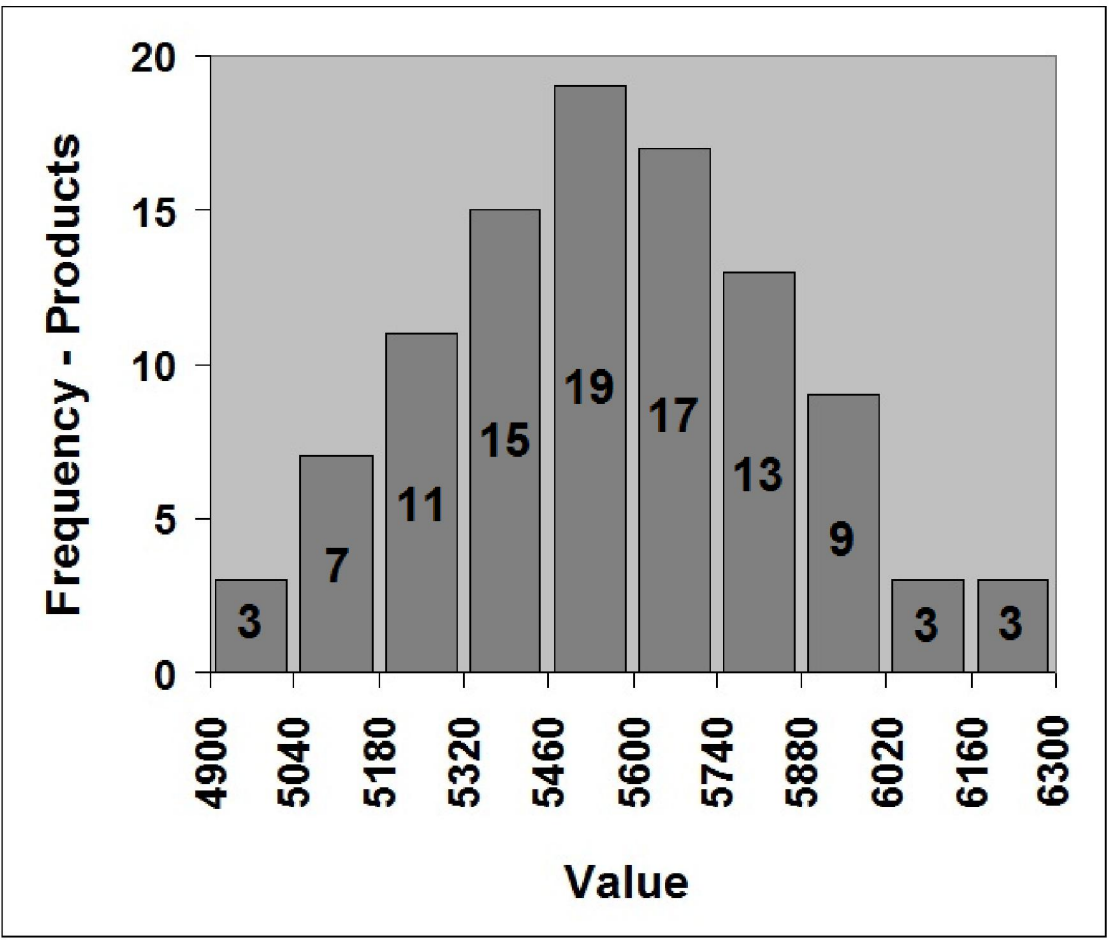

**Figure 12**: Histogram of the 70 to 79 Multiplication Table - Normal-Like Symmetry

The maximum value in the table of Figure 11 is 6241, and the minimum value is 4900; hence resultant post-multiplicative POM value for the entire table is (6241)/(4900) = (79*79)/(70*70) = (79/70)*(79/70) = 1.13*1.13 = 1.27, and this value represents just a slight increase from the 1.13 POM value of the original 70 to 79 set.

Surprisingly, this multiplication process is not skewed but rather quite symmetrical! Here the medium around the central part of the data dominates, and the histogram in Figure 12 is remarkably similar to the Normal Distribution! The only reason for this exception to the usual skewness emerging in almost all multiplication processes is the very low POM value of the 70 to 79 set of numbers. Other sets of numbers with low POM values under multiplicative processes (such as 90 to 99, or 130 to 139 for example) also show the same overall result, where most products fall in the central region of medium values, and where the histogram resembles the Normal Distribution a great deal.

# [6] Multiplicative CLT and the Lognormal Distribution

Let us provide a brief review of the Central Limit Theorem and Multiplicative Central Limit Theorem:

The **Normal distribution** is obtained in repeated additions of any random variable.
The **Lognormal distribution** is obtained in repeated multiplications of any random variable.

**Normal** $= X_1 + X_2 + X_3 + \ldots + X_N$, where $X_I$ are independent realizations from an identical random variable X. The Central Limit Theorem guarantees Normality for the sum in the limit as N gets large, and (almost) regardless of the distribution form or parameters of X.

The Central Limit Theorem requires three conditions that ostensibly should be met, namely that realizations from variable X should all be (I) independent from each other, (II) identically distributed as X, and (III) have finite variance. Yet there are quite a few generalizations weakening or even eliminating some or all of the three conditions. These generalizations ensure that CLT, and by extension Multiplicative CLT, have extremely wide applications for real life data.

The Lognormal distribution is defined as a random variable whose natural logarithm is the Normal distribution.

$$Lognormal(location, shape) = e^{Normal(location, shape)}$$
$$Normal(mean, s.d.) = LOG_e(Lognormal(mean, s.d.))$$

This definition together with the result of the CLT imply that the Lognormal distribution can be represented as a process of repeated multiplications of a random variable, namely as:

$$Lognormal = e^{[Normal]} = e^{[X_1 + X_2 + X_3 + \ldots + X_N]}$$

$$Lognormal = e^{X_1} e^{X_2} e^{X_3} \ldots e^{X_N}$$

This result, namely that repeated multiplications of a random variable is Lognormal in the limit as N gets large is called the **Multiplicative Central Limit Theorem (MCLT)**.

The shape and location parameters of the Lognormal are the standard deviation and the mean respectively of the 'generating' Normal distribution. Simulations of the Lognormal distribution for the purpose of examining its digital behavior yield an approximate logarithmic behavior whenever the shape parameter is roughly over 0.4, and regardless of the value assigned to the location parameter. Better results are obtained for higher values of the shape parameter. A near perfect logarithmic behavior is observed whenever shape parameter is approximately over 1.0, and this result is totally independent of the value assigned to the location parameter.

One plausible explanation for the prevalence of Benford's Law in the natural sciences is that such physical manifestations of the law are obtained via the cumulative effects of few or numerous multiplicative random factors, all of which leads to the Lognormal as the eventual distribution by way of the Multiplicative Central Limit Theorem. Scientists in each particular discipline have to agree on such a vista of the investigated phenomenon in focus for this to hold true, yet this research method hasn't yet been explored in the field of Benford's Law.

The crux of the matter in applying MCLT as an explanation of Benford's Law in the physical sciences is the pair of answers to the following two questions: Firstly, does Mother Nature multiply her measurements [I] sufficient number of times, or [II] merely two, three, or four times - in which case the use of MCLT may [on the face of it] seem unjustified. The second question is whether Mother Nature typically uses measurements with [I] high order of magnitude or [II] with low order of magnitude.

After some contemplations, and the creations of several hypothetical scenarios of random physical processes, together with some very limited explorations of the typical expressions in physics, engineering, chemistry and so forth, it seems safe to claim that Mother Nature is acting with considerable restraint and that she is typically reluctant to multiply more than say four or five measurements, although she often uses measurements with sufficiently high order of magnitude 'compensating' in a sense for her lack of multiplicands. Does such [partially optimistic and partially pessimistic] conclusion ruin the chance of utilizing multiplication processes as an explanation of the Benford phenomenon in the physical sciences? Certainly not! To see how Frank Benford overcomes Mother Nature's frugality in multiplying her quantities, one needs to differentiate between two very distinct goals, namely: (1) that the density of the logarithm of resultant ['multiplied'] data converges to the Normal, (2) that digit configuration converges to Benford. Mathematical empiricism strongly suggests that often convergence in the digital realm occurs faster than convergence in the stricter sense of MCLT where the goal is the emergence of the Normal shape for the log values. MCLT requires a minimum number of products of distributions in order to obtain something close enough to the Lognormal, but often in nature there are typically only 2, 3, or 4 such products, which is not sufficient for MCLT applications. Nonetheless, given that the individual variables are of sufficiently high order of magnitude, convergence in the digital sense is possible well before convergence to the Lognormal is achieved, enabling students of Benford's Law to let go of MCLT altogether.

For a decisive demonstration of the validity of the multiplicative approach as an explanation of Benford's Law in the physical sciences, let us consider Monte Carlo computer simulation results of the product of merely two Uniform distributions having high order of magnitude. In one simulation run with 35,000 realized values from the multiplicative process of the distributions Uniform(1, 650) and Uniform(1, 1200), first digits distribution came out as:

Uniform(1, 650)*Uniform(1, 1200)  {27.4, 19.7, 15.3, 11.4, 8.7, 6.5, 4.4, 3.3, 3.4}
Benford's Law First Digits Order  {30.1, 17.6, 12.5,  9.7, 7.9, 6.7, 5.8, 5.1, 4.6}

The relatively low 28.3 SSD value signifies that results are quite close to Benford. The Uniform(1, 650) has the high POM value of 650/1 = 650, and the Uniform(1, 1200) has the high POM value of 1200/1 = 1200, and such high orders of magnitude guarantee that their product is very close to the Benford digital configuration.

Yet, multiplications of two Uniform distributions having exceptionally low order of magnitude does not yield any convergence to Benford whatsoever. In one simulation run with 35,000 realized values from the multiplicative process of the distributions Uniform(5, 11) and Uniform(8, 31), first digits distribution came out as:

Uniform(5, 11)*Uniform(8, 31)  {52.7, 23.0, 1.9, 0.8, 2.5, 3.4, 4.5, 5.7, 5.6}
Benford's Law First Digits Order  {30.1, 17.6, 2.5, 9.7, 7.9, 6.7, 5.8, 5.1, 4.6}

The very high 790.3 SSD value signifies that results are decisively non-Benford. The Uniform(5, 11) has the very low POM value of 11/5 = 2.2, and the Uniform(8, 31) has the very low POM value of 31/8 = 3.9, and such exceedingly low orders of magnitude preclude convergence to the Benford digital configuration.

The first result regarding two Uniform distributions with high orders of magnitude is a decisive indication of how closely results can get to Benford in merely a single multiplicative process having only two multiplicands. In other words, that Benford behavior can be achieved quite easily and rapidly, long before resultant data becomes Lognormal via MCLT.

Conceptually, as a general rule, Benford digital configuration is expected to be found in the natural sciences whenever **'the randoms are multiplied'**. More specifically, the multiplicative process needs to be carefully scrutinized to determine whether there are sufficiently many multiplicands, or at least sufficiently high order of magnitude within the constituent distributions. The crucial factor here is the resultant order of magnitude of the entire multiplicative process, which depends on the individual orders of magnitude of the distributions being multiplied, as well as on the number of distributions involved. The higher the number of distributions being multiplied and the higher the orders of magnitude of the various individual distributions – the larger is resultant order of magnitude of the multiplicative process itself and the closer is resultant digital configuration to Benford.

# [7]  The Effects of Random Multiplication Processes on POM

Let us examine theoretically how multiplication processes of random variables affect POM. Given PDF(x), namely Probability Density Function of variable x with $MIN_X$ and $MAX_X$ endpoints; and given PDF(y), namely Probability Density Function of variable y with $MIN_Y$ and $MAX_Y$ endpoints; then   $POM_X = MAX_X / MIN_X$   and   $POM_Y = MAX_Y / MIN_Y$.

For the random multiplication process X*Y:
The best scenario is multiplying $MAX_X*MAX_Y$
The worst scenario is multiplying $MIN_X*MIN_Y$
These two extreme scenarios represent the maximum and minimum values for the multiplication process itself.

Hence, for the random multiplicative process PDF(x)*PDF(y):
$POM_{X*Y} = [ MAX_X*MAX_Y ]/[ MIN_X*MIN_Y ] = [ MAX_X / MIN_X ]*[ MAX_Y / MIN_Y ] = POM_X*POM_Y$. Since by definition POM  > 1 for each variable that is not a constant, it follows that there is a definite increase in resultant POM due to the multiplicative act.

In reality, for any finite number of computer simulated multiplications of two random variables, such as in say thousands or tens of thousands realizations, it is rare to simultaneously obtain the highest values of the two in one single realization; and it is also rare that the two lowest values are gotten simultaneously. The maximum multiplied value in the entire simulation process is a bit lower than the theoretical; and the minimum multiplied value in the entire simulation process is a bit higher than the theoretical; all of which implies that resultant POM of multiplied variables is a bit lower in actual simulations then the theoretical expression $POM_X*POM_Y$.

Yet, the basic argument presented here is certainly valid, and actual resultant $POM_{X*Y}$ strongly increases (from the $POM_X$ level and from the $POM_Y$ level) due to the multiplicative process even if it always hovers somewhere below the theoretical $POM_X*POM_Y$ expectation.

Needless to say, this result could easily extend to any number of random variables, not merely two, and where resultant theoretical POM is the product of all the individual POMs, namely $POM_{PRODUCT} = POM_X*POM_Y*POM_Z$ and so forth; or in concise mathematical notations as:

$$POM_{PRODUCT} = \prod POM_J$$

For the 10 by 10 multiplication table from our elementary school days, POM of each multiplicand {1, 2, 3, 4, 5, 6, 7, 8, 9, 10} is 10/1 or simply 10, while POM of the entire table is 100/1 or simply 100. This result nicely complies with the above derived expression of $POM_{X*Y} = POM_X*POM_Y = 10*10 = 100$.

For statisticians who are more accustomed to the logarithm definition of OOM as discussed in chapter IV of Part 1, the net effect of multiplication processes is simply the summing up of the two distinct orders of magnitude as in $OOM_{X*Y} = OOM_X + OOM_Y$, leading to the same conclusion, namely that all multiplication processes yield increased OOM and variability.

# [8] Scrutinizing Random Multiplication Processes

Let us examine multiplication processes of Uniform distributions via Monte Carlo computer simulations in order to learn how they gradually build up sufficiently large cumulative POM and increase the system's skewness; thereby converging to Benford.

*[Note: all simulations are with 35,000 realized values.]*

---

Uniform(3, 40) POM = 13
1st significant digits: {26.7, 27.4, 29.7, 2.6, 2.9, 2.7, 2.5, 2.8, 2.7} SSD = 514

---

Uniform(3, 40)*Uniform(2, 33) POM = 212
1st significant digits: {23.5, 15.9, 13.1, 11.2, 9.3, 8.5, 7.2, 6.1, 5.1} SSD = 58

---

Uniform(3, 40)*Uniform(2, 33)*Uniform(7, 41) POM = 1076
1st significant digits: {32.2, 18.1, 11.3, 8.4, 7.2, 6.6, 5.9, 5.5, 4.9} SSD = 9

---

Uniform(3, 40)*Uniform(2, 33)*Uniform(7, 41)*Uniform(1, 29) POM = 5386
1st significant digits: {30.3, 18.0, 12.9, 9.5, 7.7, 6.5, 5.6, 4.9, 4.6} SSD = 1

---

| Random Process | POM Variability/Range | SSD Deviation from Benford |
|---|---|---|
| U | 13 | 514 |
| U*U | 212 | 58 |
| U*U*U | 1076 | 9 |
| U*U*U*U | 5368 | 1 |

**Figure 13**: Increase in POM Variability Measure Induces Closeness to Benford

The table in Figure 13 demonstrates how larger POM value goes hand in hand with *[induces rather]* closeness to Benford, with the lowering of SSD at each stage of the multiplicative process. The sequence of empirical POM values in Figure 13 can be easily checked against the theoretical expression $POM_{PRODUCT} = \prod POM_J$ since the max and the min of Uniforms can be readily obtained as in Uniform(min, max). Referring to the variables Uniform(3, 40), Uniform(2, 33), Uniform(7, 41), Uniform(1, 29); the four individual POM values in the same order as above are {40/3, 33/2, 41/7, 29/1}, namely {13.3, 16.5, 5.9, 29.0}. Successive multiplications yield {13.3, 13.3*16.5, 13.3*16.5*5.9, 13.3*16.5*5.9*29.0}, namely the POM sequence of **{13, 220, 1289, 37369}**. But this theoretical sequence of sharp rise in POM is not what seems in Figure 13. As mentioned above, it is rare to simultaneously obtain the highest/lowest values in any one single realization in actual simulations; hence we empirically encounter the somewhat less dramatic rise in POM of the sequence **{13, 212, 1076, 5368}** as seems in the empirical results of Figure 13.

The histogram in Figure 14 of the log values of this multiplicative process demonstrates in general what could occur at times in random multiplication processes, namely that digits converge to Benford well before any significant achievement of MCLT in terms of endowing log density the shape of the Normal. Here, log of U(3, 40)*U(2, 33)*U(7, 41)*U(1, 29) is a bit asymmetrical with a longer tail to the left, not resembling sufficiently the shape of the Normal distribution, yet digits here are nearly perfectly logarithmic, with an exceedingly low SSD value of 1. Such low SSD value is rarely found in real-life random data, and it indicates nearly perfect agreement with Benford's Law.

*[CLT and MCLT strictly requires 'independent and identically distributed variables', while here we consider the product of four non-identical distributions. Nonetheless, the many extensions and generalizations of CLT and MCLT allow us to consider such non-identical distributions as well, so the general principle holds.]*

The histogram in Figure 15 of the actual simulation values (not log-transformed) depicts the severe quantitative skewness of the resultant distribution. It is exactly this quantitative skewness which drives digits into their Benford digital skewness as a consequence.

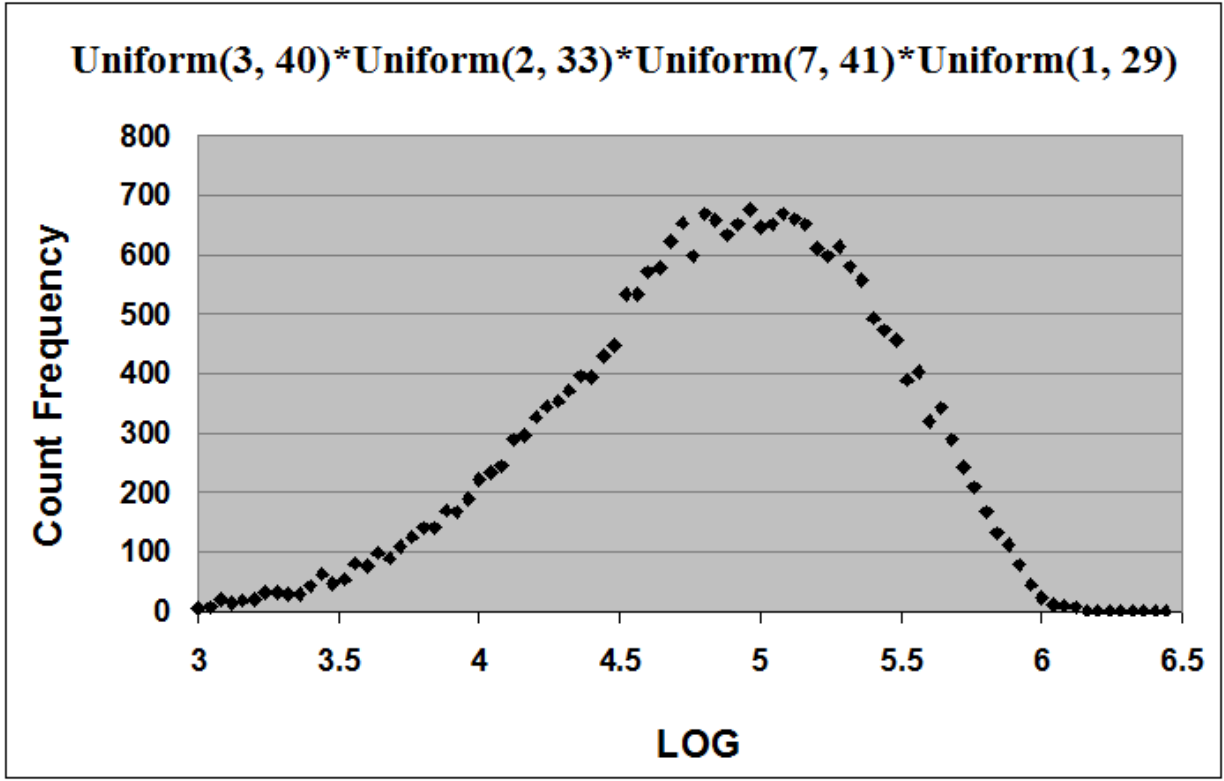

**Figure 14**: Log Histogram of U(3, 40)*U(2, 33)*U(7, 41)*U(1, 29) is not Quite Normal

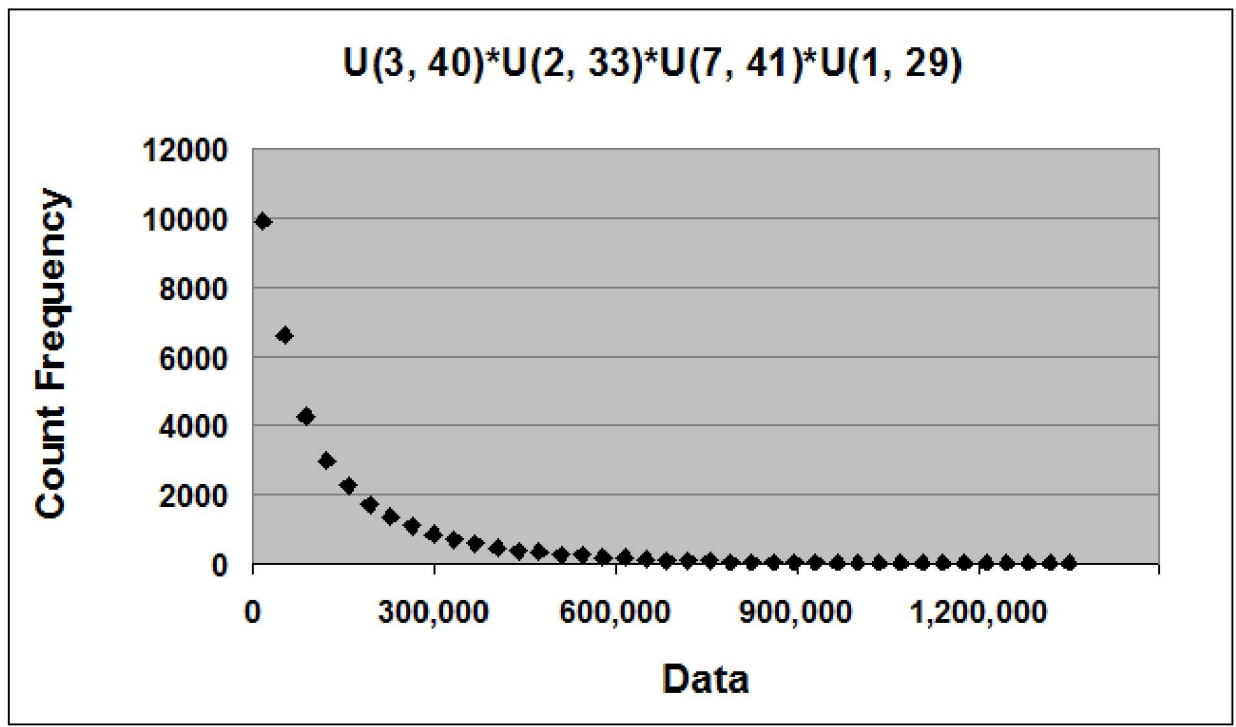

**Figure 15**: Histogram of U(3, 40)*U(2, 33)*U(7, 41)*U(1, 29) is Highly Skewed and Benford

When these Uniform distributions themselves are of **very high POM,** merely a single product of two such high-variability Uniforms yields digital configuration quite close to Benford, all in one fell swoop. As an example, Monte Carlo computer simulations are performed for the product of Uniform(1, 60777333) and Uniform(1, 30222888). Theoretical POM for these two Uniforms are 60777333 and 30222888 respectively. Theoretical POM for their product is 60777333*30222888 = 1836866528197700, although actual simulations of the multiplication process should yield much lower POM value.

Five simulations are performed for Uniform(1, 60777333)*Uniform(1, 30222888); with 35,000 runs each; the results are:

First Digits: {29.7, 13.6, 11.5, 10.3, 8.7, 7.6, 6.8, 6.5, 5.3}   SSD = 22.5   POM =     661343
First Digits: {29.0, 13.5, 11.6, 10.3, 9.0, 8.1, 7.0, 5.9, 5.6}   SSD = 25.9   POM =     455568
First Digits: {29.4, 14.0, 11.5, 10.1, 8.9, 7.6, 7.0, 6.1, 5.4}   SSD = 19.5   POM = 26872700
First Digits: {29.2, 13.3, 11.8, 10.1, 9.1, 7.8, 6.9, 6.2, 5.7}   SSD = 26.0   POM =     114111
First Digits: {29.4, 13.6, 11.5, 10.3, 8.7, 7.7, 6.8, 6.3, 5.5}   SSD = 22.8   POM =   3550739

Yet, MCLT does not even begin to be significantly effective for this very short multiplicative process! The histogram in Figure 16 of the log values of this process demonstrates once again that in multiplication processes digits could converge to Benford well before any significant achievement of MCLT is obtained in terms of endowing log density the shape of the Normal distribution (or equivalently that resultant distribution is Lognormal). Here, log histogram of Uniform(1, 60777333)*Uniform(1, 30222888) does not even begin to resemble the Normal, and yet digits here are quite close to the logarithmic, with the relatively low SSD value of about 22!

The histogram in Figure 17 of the actual simulation values (not log-transformed) depicts the quantitative skewness of the resultant distribution. It is exactly this quantitative skewness which drives digits into their Benford digital skewness as a consequence.

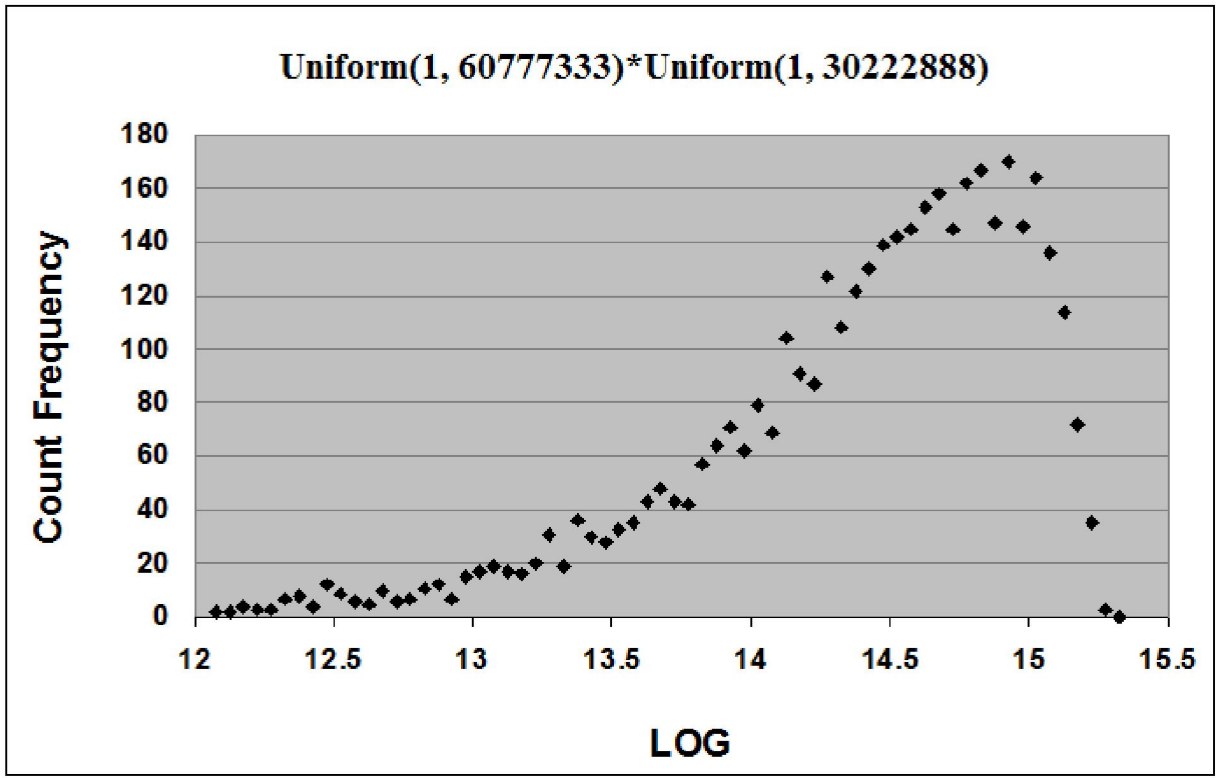

**Figure 16**: Log of U(1, 60777333)*U(1, 30222888) is Totally not Normal but Data is Benford

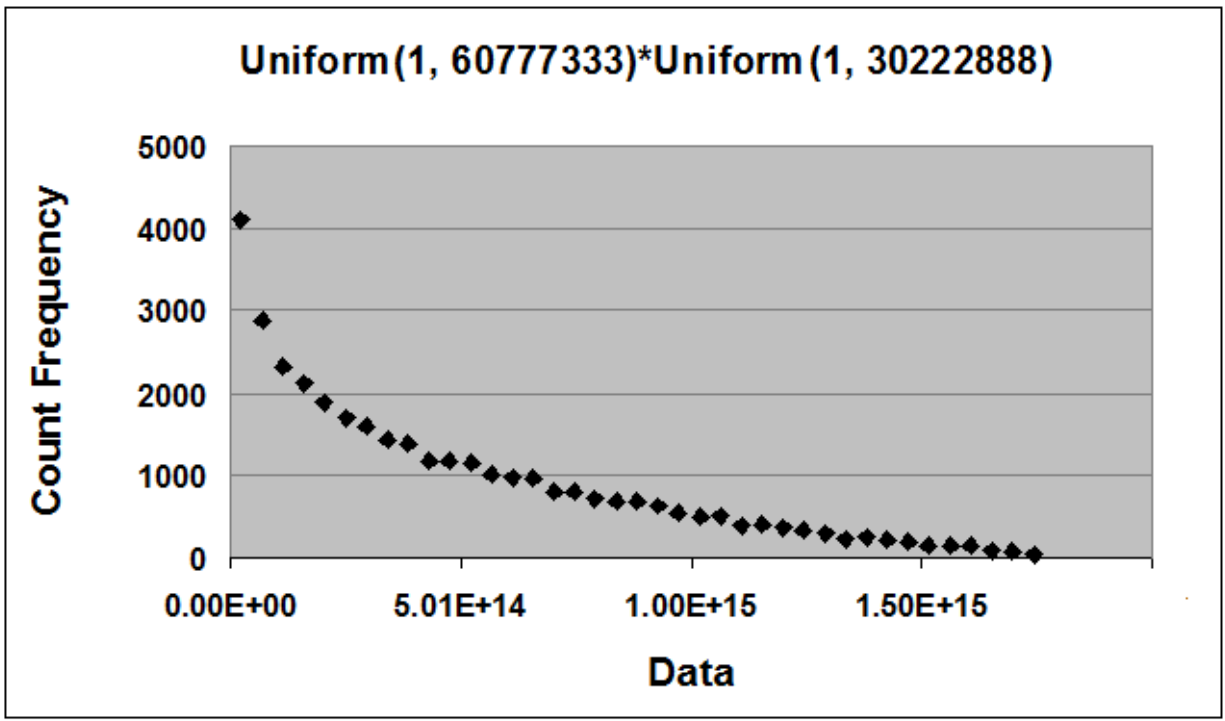

**Figure 17**: Histogram of U(1, 60777333)*U(1, 30222888) is Quite Skewed and Benford

When these Uniform distributions themselves are of **very low POM,** not even the product of several such low-variability Uniforms is capable of yielding digital configuration close enough to Benford, unless a truly great many such Uniforms are multiplied. As an example, Monte Carlo simulations are performed for the product of 6 Uniforms with very low POM as follow: Unif(4, 7)*Unif(8, 11)*Unif(5, 7)*Unif(12, 16)*Unif(237, 549)*Unif(17, 25).
Theoretical POMs for these six Uniforms are 7/4, 11/8, 7/5, 16/12, 549/237, 25/17, namely 1.8, 1.4, 1.4, 1.3, 2.3, 1.5, respectively. Theoretical POM for their product is
1.8*1.4*1.4*1.3*2.3*1.5 = 15.3, although actual simulations of the multiplication process should yield much lower POM value.

Five simulations for these six Uniforms are conducted; with 35,000 runs each; the results are:

First Digits: {6.7, 28.5, 30.8, 19.7, 9.5, 3.5, 1.0, 0.2, 0.0}   SSD = 1181.1   POM = 11
First Digits: {6.6, 29.0, 30.9, 19.6, 9.5, 3.3, 0.9, 0.2, 0.0}   SSD = 1197.5   POM =  9
First Digits: {6.7, 28.8, 30.6, 19.8, 9.4, 3.5, 1.0, 0.2, 0.0}   SSD = 1183.3   POM = 10
First Digits: {6.5, 28.9, 30.7, 19.9, 9.3, 3.4, 1.0, 0.2, 0.0}   SSD = 1202.1   POM = 10
First Digits: {6.5, 28.8, 31.2, 19.6, 9.3, 3.3, 1.0, 0.2, 0.0}   SSD = 1213.6   POM =  9

Here, MCLT is properly and nicely applied for this longer process with sufficient number of multiplications (almost), and the histogram of the log is very much Normal-like in appearance. The resultant data set itself is nearly Lognormal, but with a low shape parameter – precluding Benford behavior. Figure 18 of the histogram of the log of this process demonstrates what could happen at times in multiplication processes of random variables having very low POM values. Here digits do not converge to Benford in the least, while the significant achievement of MCLT is obvious and visible, endowing log density the shape of the Normal (almost).

The histogram in Figure 19 of the actual simulation values (not log-transformed) depicts the approximate quantitative symmetry of resultant distribution. The non-skewed quantitative nature of resultant data precludes any resemblance to the Benford skewed digital configuration.

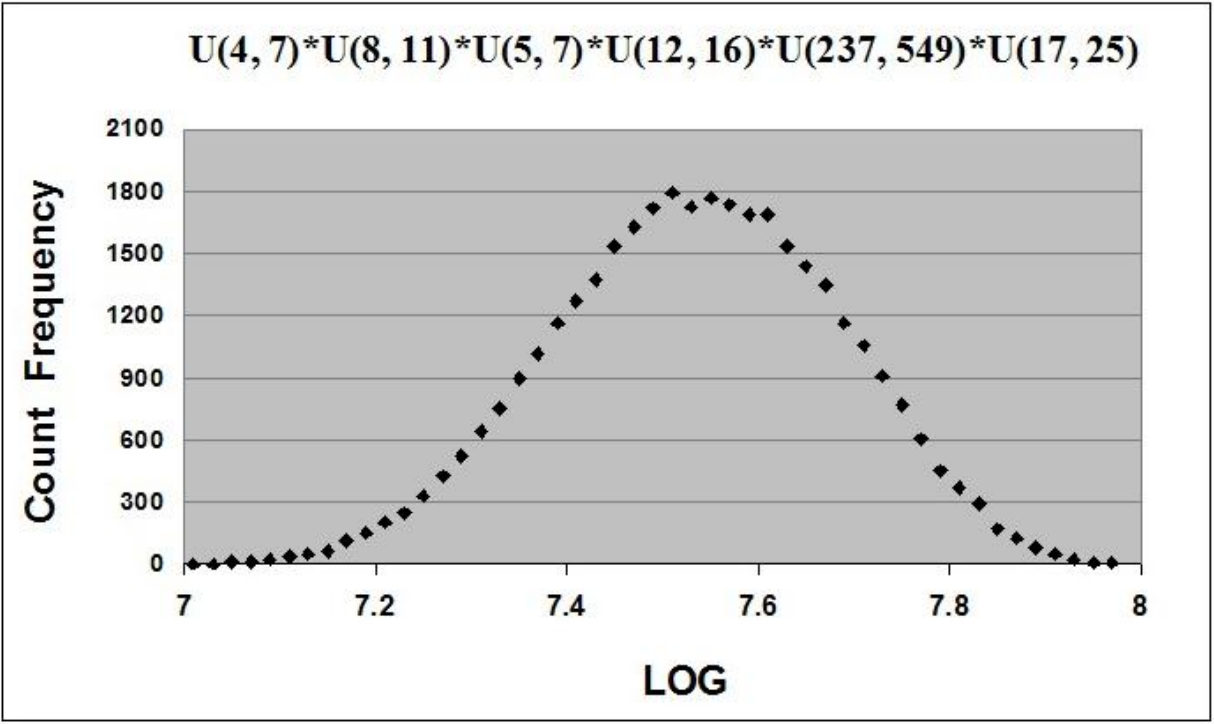

**Figure 18**: Log of Six Multiplied Uniforms with Low POM is Normal but Data is not Benford

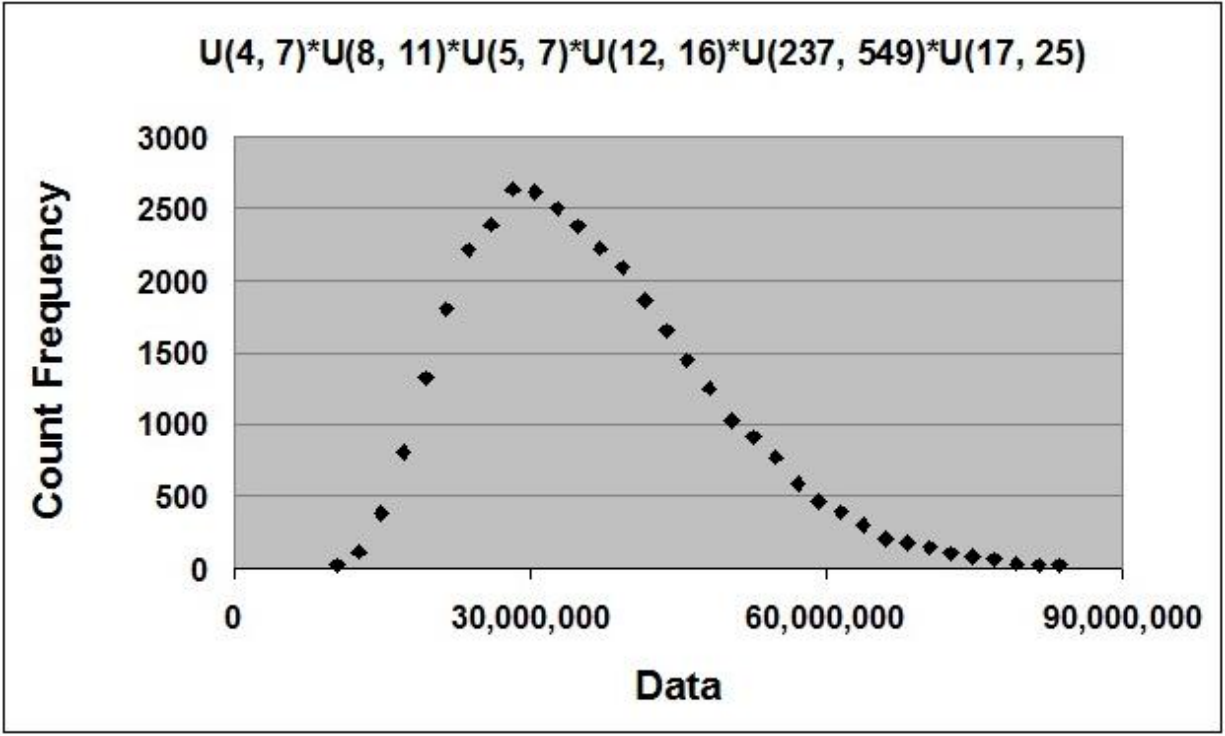

**Figure 19**: Histogram of Six Multiplied Uniforms with Low POM is Nearly Symmetrical

The essential lessons that can be learnt from the three examples above can be stated succinctly:

**1)** MCLT does not guarantee Benford behavior in multiplication processes!
**2)** MCLT is not necessary for Benford behavior in multiplication processes!
**3)** Order of Magnitude is the crucial factor in determining Benford behavior in multiplication processes!
**4)** The confluence of two factors leads to high resultant POM and thus to Benford behavior and skewness in multiplication processes: (a) POMs of the individual variables, (b) the number of variables that are being multiplied.

POM is a greedy creature when it comes to multiplication processes, it is never satisfied with what it has, hence the more variables we multiply the larger POM gets with each added multiplicand, and indefinitely so. MCLT on the other hand is modest; it has a limited and well-defined goal, namely, Normality of the density of the logarithm. Benford is also very modest; it has a limited and well-defined goal, namely LOG(1 + 1/d). Once log Normality is approximately achieved for MCLT - it is maintained, and throwing in more variables helps only slightly in terms of perfection of convergence. Once LOG(1 + 1/d) is approximately achieved for Benford - it is maintained, and throwing in more variables helps only slightly in terms of perfection of convergence. The grand generic scheme of adding more and more multiplicands to the process (namely multiplying more and more variables one by one in stages) could lead to either MCLT achieving its goal first with log being nearly Normal, while digits are not yet Benford [*the case of low POM values of many multiplicands*], or it could lead to Benford achieving its goal first and digits are nearly LOG(1 + 1/d), while log curve is not yet Normal [*the case of high POM values of few multiplicands*]. But no matter who is first, the grand multiplicative scheme eventually leads to both, to MCLT as well as to Benford, both ultimately achieving their respective goals with certainty. With delicate calibration one could perhaps choose variables/multiplicands that let both MCLT and Benford achieve their respective goals simultaneously at approximately the same stage; namely after approximately the same number of multiplications – but this issue is neither significant nor interesting.

Hence, in light of the above discussion, instead of desperately seeking justifications and validations via the *[often unrealistic or irrelevant]* Multiplicative CLT, students of Benford's Law should learn to rely instead on *[the very realistic, relevant, and easily available]* Related Log Conjecture [chapter V of Part 1 and Kossovsky (2014) chapter 63] and focus on resultant order of magnitude which is the range or width of the histogram of the logarithm of the data. Examinations of several log histograms of multiplied distributions (with sufficient POM) limited to only 2, 3, or 4 such products, show that while they are not really Normal-like, nonetheless the requirements of Related Log Conjecture are almost all nicely met, and thus digits are nearly Benford.

A decisive argument supporting the above discussion can be found in the consideration of log of data in the form of a triangle which yield the Benford configuration as close to perfection as can be measured with any finite set of simulated values for $10^{\text{Triangular}}$ whenever log-range is say over 4 approximately [see Kossovsky (2014) chapters 64 and 66 for details]. Moreover, Lawrence Leemis (2000) rigorously proved that a symmetrical triangle positioned onto integral log points is exactly logarithmic. Here the Benford configuration is obtained so perfectly, and without log density being Normal-like in any way. To emphasis this point once again: there is no need whatsoever for log of data to be Normal in order to obtain an exact or near-perfect Benford behavior! Benford's Law and Normality of log are not two sides of the same coin! True, numerous *[Lognormal-like]* physical and scientific data sets which are Benford come with log density that is almost Normal, or at least resembles the Normal a great deal, but this is not what Benford's Law is all about. The law is much more general than this narrow manifestation of it.

Two manifestations of some very short multiplication processes having only 3 multiplicands and resulting in near Benford configuration demonstrate that such situations can occur in scientific and physical data quite frequently. These two scenarios are discussed in Kossovsky (2014) chapter 90 page 393.

Case I: The final speed of numerous particles of random mass M, initially at rest, driven to accelerate linearly via random force F, applied for random length of time T. In other words, each particle comes with (random) unique M, F, and T combination. The equations in Physics describing the motion are: (i) $V_{FINAL} = V_{INITIAL} + A*T$, or since they start at rest simply $V_{FINAL} = A*T$, (ii) $F = M*A$. Thus $V_{FINAL} = (Force*Time)/Mass$. Monte Carlo computer simulations show strong logarithmic behavior here whenever M, F, and T are randomly chosen from Uniform distributions of the form U(0, b) *[intentionally choosing parameter a equals to 0 to ensure high order of magnitude]*. Moreover, even with only 2 random measurements, and with the 3rd measurement being fixed as a constant, we get quite close to the logarithmic [say mass M and force F are random, while time T is fixed - being a constant for all particles.]

Case II: The final position of numerous particles of random mass M, thrown linearly from the same location at random initial speed $V_I$, under constant decelerating random force F (say frictional) until each one comes to rest. In other words, each particle comes with (random) unique M, F, and $V_I$ combination. The equations in Physics describing the motion are:
(i) $V_{FINAL} = V_I + A*T$, which leads to $T_{REST} = V_I/A$ as the time it takes to achieve rest,
(ii) $F = M*A$, (iii) $Displacement = V_I*T - (1/2)*A*T^2$. Hence
$Displacement = V_I*(V_I/A) - (1/2)*A*(V_I/A)^2$, namely $Displacement = (V_I)^2/(2*(F/M))$. Monte Carlo simulations again show strong logarithmic behavior here as well whenever M, F, and $V_I$ are randomly chosen from Uniform(0, b) *[intentionally choosing parameter a equals to 0 to ensure high order of magnitude]*. Moreover, even with only 2 random measurements, and with the 3rd measurement being fixed as a constant, we get quite close to the logarithmic [say mass M and initial speed $V_I$ are random, while frictional force F is fixed - being a constant for all particles.]

Schematically Case I is of the form: Simulation = (U1*U2)/U3, while Case II is of the form: Simulation = (U1*U1)/(2*(U2/U3)), and so neither can really apply MCLT due to severe lack of multiplicands; nevertheless due to high order of magnitude of multiplied variables digits here are quite close to the logarithmic.

As an additional example of the rapidity and the decisiveness with which multiplication processes could lead to Benford, Monte Carlo computer simulation results of the product of just two Normal(m, sd) distributions would be shown. The Normal is considered to be 'anti-logarithmic' in and of itself, yet as multiplicands it can readily become logarithmic.

Normal(2, 9)*Normal(5, 13)  gave  {31.7, 17.6, 11.3,  9.7, 7.8, 6.0, 5.8, 5.4, 4.9}
Normal(4, 7)*Normal(2, 3)   gave  {28.7, 18.5, 13.1, 10.2, 7.7, 6.6, 5.9, 4.9, 4.3}
Normal(2, 4)*Normal(5, 3)   gave  {29.0, 19.3, 13.3, 10.6, 8.1, 6.4, 4.7, 4.9, 3.7}
-------------------------------------------------------------------------------------------------
Benford's Law First Digits         {30.1, 17.6, 12.5,  9.7, 7.9, 6.7, 5.8, 5.1, 4.6}

Admittedly, the choices of parameters above were somewhat intentional, ensuring that the range is of large order of magnitude *[accomplished by crossing the origin, which ensures that the data draws plenty from the interval (0, 1)]*, but the general principle here holds nonetheless, and certainly it could be argued that in nature typical Normals are such, often occurring with high order of magnitude.

Monte Carlo computer simulation results of the product of just two Exponential distributions, a density form already somewhat close to Benford in and of itself, yield much superior results:

Exponential(4)*Exponential(11.0) gave  {30.1, 17.6, 12.8, 9.9, 7.3, 7.1, 6.0, 4.7, 4.5}
Exponential(5)*Exponential(0.07) gave  {30.2, 17.5, 12.7, 9.3, 9.4, 7.3, 5.8, 3.8, 3.9}
Exponential(13)*Exponential(0.2) gave  {30.4, 17.1, 12.9, 9.9, 7.5, 6.9, 5.7, 5.3, 4.3}
-------------------------------------------------------------------------------------------------
Benford's Law First Digits              {30.1, 17.6, 12.5, 9.7, 7.9, 6.7, 5.8, 5.1, 4.6}

# [9] The Conjecture on the Benfordized Nature of Multiplication

In Kossovsky (2014) chapter 58, two important results by Hamming (1970) are mentioned relating to products and divisions of several data sets or distributions. Let X be a continuous random variable with an exact logarithmic behavior. Let Y be any other continuous random variable, logarithmic or non-logarithmic. Then the product X*Y, and the ratios X/Y, Y/X, all satisfy Benford's Law as well. The ramification of these two remarkable properties is profound, since it applies to any distribution form Y! These two properties may be thought of as 'propagators' of the logarithmic distribution, guaranteeing to spread it around whenever a logarithmic data set gets 'multiplicatively connected' with any other data type, with the result that in turns it 'infects' other data sets with the logarithmic configuration, and so on. Consequently, in the context of scientific and physical data sets, even 2 or 3 products are sufficient to obtain a perfect Benford behavior given that one measurement is Benford in its own right! Moreover, one can conjecture with total confidence that partial Benfordness in one measurement endows *[at least]* that same degree of Benfordness *[and actually somewhat higher degree]* for the product with any other measurement. This extrapolation for products of random variables is in the same spirit as the extrapolation of the 2nd chain conjecture in Kossovsky (2014) chapter 102 page 460, where it is claimed that any single sequence within a chain of distribution contributes to increase in Benfordness, or that at least there is no backtracking away from it.

If one goes along with SSD as a measure of 'distance' from the logarithmic, then this conjecture can be stated succinctly and formally as:

[ SSD of X*Y ] $\leq$ [ SSD of X ]
                and
[ SSD of X*Y ] $\leq$ [ SSD of Y ]

for any X, Y random variables or data sets of positive numbers.

The mixed inequality/equality sign should be substituted with inequality sign exclusively whenever X or Y are not exactly (perfectly) Benford, signifying an improvement in Benfordness for any product of variables whatsoever. The equality sign is valid for cases where X or Y are exactly Benford to begin with prior to any multiplication, so that SSD is 0 on both sides.

# [10]  Random Additions of Lognormal Distributions

In sharp contrast to multiplication processes, random additions of random variables affect digital configuration in the opposite way, away from Benford. Random addition processes induce symmetry in resultant histogram, gradually eliminating skewness altogether, and in the limit the emerging shape of resultant histogram is the bell curve of the Normal Distribution.

In order to emphasize that **addition is highly detrimental to Benford behavior**, eight highly logarithmic Lognormal distributions with shape parameter well over 1 are added via Monte Carlo computer simulations, one by one, resulting in significant deviations and retreat from Benfordness - as predicated by the Central Limit Theorem. *[Note: all simulations of Lognormals and their sums are with 35,000 realized valued.]*

Lognormal( shape = 1.5,  location = 3.8)
First Digits: {29.7, 17.6, 12.8, 9.6, 8.0, 6.8, 5.8, 5.1, 4.6}     SSD = 0.2

Lognormal( shape = 1.3,  location = 4.0)
First Digits: {29.9, 17.8, 12.4, 9.6, 7.8, 7.0, 5.8, 5.1, 4.6}     SSD = 0.2

Lognormal( shape = 1.6,  location = 4.3)
First Digits: {30.5, 17.5, 12.6, 9.6, 7.8, 6.4, 5.7, 5.1, 4.8}     SSD = 0.3

Lognormal( shape = 1.2,  location = 4.2)
First Digits: {30.3, 17.5, 12.6, 9.9, 7.9, 6.7, 5.7, 4.9, 4.6}     SSD = 0.1

Lognormal( shape = 1.4,  location = 4.1)
First Digits: {30.2, 17.4, 12.9, 9.7, 8.0, 6.6, 5.8, 4.9, 4.5}     SSD = 0.3

Lognormal( shape = 1.1,  location = 4.4)
First Digits: {30.4, 17.0, 12.1, 9.7, 8.2, 6.9, 5.9, 5.2, 4.6}     SSD = 0.7

Lognormal( shape = 1.3,  location = 4.3)
First Digits: {30.3, 17.4, 12.4, 9.9, 7.8, 6.7, 5.6, 5.1, 4.7}     SSD = 0.2

Lognormal( shape = 1.5,  location = 4.5)
First Digits: {30.3, 17.4, 12.4, 9.9, 7.8, 6.7, 5.6, 5.1, 4.7}     SSD = 0.2

-----------------------------------------------------------------------------------------

**Benford 1st:** {30.1, 17.6, 12.5,  9.7, 7.9, 6.7, 5.8, 5.1, 4.6}     SSD = 0

LogN(1.5, 3.8) + LogN(1.3, 4.0)
{31.6, 17.8, 11.7, 9.1, 7.4, 6.4, 5.9, 5.2, 5.0}
SSD = **3.7**

LogN(1.5, 3.8) + LogN(1.3, 4.0) + LogN(1.6, 4.3)
{30.0, 19.3, 13.2, 9.4, 7.6, 6.3, 5.2, 4.9, 4.2}
SSD = **4.2**

LogN(1.5, 3.8) + LogN(1.3, 4.0) + LogN(1.6, 4.3) + LogN(1.2, 4.2)
{25.5, 19.1, 14.9, 11.2, 8.4, 6.9, 5.6, 4.5, 3.9}
SSD = **32.6**

LN(1.5, 3.8) + LN(1.3, 4.0) + LN(1.6, 4.3) + LN(1.2, 4.2) + LN(1.4, 4.1)
{23.5, 15.9, 14.5, 12.0, 9.7, 8.0, 6.7, 5.3, 4.5}
SSD = **60.1**

LN(1.5, 3.8) + LN(1.3, 4.0) + LN(1.6, 4.3) + LN(1.2, 4.2) + LN(1.4, 4.1) + LN(1.1, 4.4)
{25.0, 11.7, 11.7, 11.6, 10.7, 9.2, 8.1, 6.5, 5.6}
SSD = **87.2**

LN(1.5, 3.8)+LN(1.3, 4.0)+LN(1.6, 4.3)+LN(1.2, 4.2)+LN(1.4, 4.1)+LN(1.1, 4.4)+LN(1.3, 4.3)
{31.1, 9.7, 8.4, 9.3, 9.9, 9.0, 8.6, 7.4, 6.7}
SSD = **107.7**

L(1.5, 3.8)+L(1.3, 4.0)+L(1.6, 4.3)+L(1.2, 4.2)+L(1.4, 4.1)+L(1.1, 4.4)+L(1.3, 4.3)+L(1.5, 4.5)
{38.2, 11.6, 6.4, 6.5, 7.4, 7.8, 7.8, 7.3, 7.0}
SSD = **166.7**

The additions of these eight Lognormal distributions constitute in essence a tug of war between additions and multiplications, a war decisively won by additions, overcoming multiplications, thus disobeying the law of Benford. This is so since each Lognormal distribution may be represented as a random multiplicative process. For example, assuming that the product of only 6 Uniforms is a sufficient approximation of a Lognormal distribution, then the addition of the 3 Lognormals above LogN(1.5, 3.8) + LogN(1.3, 4.0) + LogN(1.6, 4.3) could be represented by the expression U*U*U*U*U*U + U*U*U*U*U*U + U*U*U*U*U*U where the tension between multiplication and addition is clearly demonstrated.

Figure 20 vividly illustrates this tension that exists between additions and multiplications in the context of Benford's Law. Assuming for the sake brevity and lack of space that the product of only 3 Uniforms is a sufficient approximation of a Lognormal distribution (it is not!), namely that Lognormal ≈ U*U*U, then the above arrangement of successive additions of Lognormals, gradually leading to significant deviations from Benford, can be viewed as a gradual victory by additions over multiplications due to the fact that at each successive stage these victorious additions are a notch more overwhelming or numerous than the vanquished multiplications. The damage done by additions here is mostly due to the gradual transformation of the shape of the resultant curve from skewness to symmetry, and not so much from the gradual decrease in POM. Indeed, the final value of 783 POM is still large enough and almost Benford in potential.

| The Random Process | POM | SSD |
|---|---|---|
| UUU | 1133924 | 0.2 |
| UUU+UUU | 37836 | 3.7 |
| UUU+UUU+UUU | 14914 | 4.2 |
| UUU+UUU+UUU+UUU | 4747 | 32.6 |
| UUU+UUU+UUU+UUU+UUU | 2272 | 60.1 |
| UUU+UUU+UUU+UUU+UUU+UUU | 1416 | 87.2 |
| UUU+UUU+UUU+UUU+UUU+UUU+UUU | 979 | 107.7 |
| UUU+UUU+UUU+UUU+UUU+UUU+UUU+UUU | 783 | 166.7 |

**Figure 20**: Tug of War between Additions and Multiplications

# [11]  The Achilles' Heel of the Central Limit Theorem

Yet, the full description of such arithmetical tugs of war is a bit more complex than the apparently simple narrative outlined above in Figure 20. First of all, one should bear in mind that CLT works decisively and fast whenever added variables are non-skewed (symmetric). When added variables are skewed, more 'work' (in terms of adding very many of them) is required in order to shape their sum into that symmetric Normal curve from the original non-symmetric configuration. The **CLT's Achilles' heel** - in terms of its rate of convergence to the Normal - is the adverse possibility that added variables are highly skewed and come with very high order of magnitude, constituting a very bad combination for the CLT and a challenge that needs to be overcome. Except for Uniforms, Normals, and other symmetric distributions which converge to the Normal quite fast after very few additions regardless of the value of order of magnitude of added variables, all other asymmetrical (skewed) distributions show a distinct rate of convergence depending on the value of their order of magnitude. For skewed variables, whenever order of magnitude is of very high value, CLT can manifest itself with difficulties, and very slowly, only after a truly large number of additions of the random variables. On the other hand, when skewed variables are of very low order of magnitude, CLT achieves near Normality quite quickly after only very few additions.

As an example of the challenge CLT faces when added variables are highly skewed and of very high order of magnitude; 8 identical Lognormals with very high shape parameter of 2.2 and location parameter 5 are randomly added. The shape parameter of the Lognormal determines its order of magnitude, and for shape value of 2.2 order of magnitude is almost 6 and POM is approximately 1,000,000. These 8 Lognormals are thus with larger POM values than the 8 Lognormals of Figure 20 which were with lower shape parameter value of around 1.36 and thus of lower POM. The shape of the resultant data set after the 8 additions of these high POM Lognormals in 35,000 simulation runs shows a decisively non-Normal and highly skewed curve much like the original Lognormal; so that not much has been achieved after 8 additions, and which is in sharp contrast to what is shown in Figure 20. Where is CLT? What happened to the mathematicians' solemn pledge of seeing the Normal curve after several random additions? Well, here, when dealing with highly skewed variables with very high POM values, a lot of patience and very strong faith in CLT is necessary. The Average of resultant added data of this simulation is 13325, being much greater than the Median value of 4055, and which is a strong indication of positive skewness. Leading digits of the 8 added Lognormals are still very much Benford-like coming at {29.9, 17.8, 12.5, 9.7, 7.9, 6.7, 5.7, 5.2, 4.7} with an extremely low SSD value of 0.1. It would take many more additions of these highly skewed Lognormals with high POM to arrive approximately at the Normal, and surely our (rational) faith in CLT would eventually be rewarded.

On the other hand, Lognormals with much lower shape parameter (and the implied lower POM value) converge to the Normal quite rapidly in spite of their skewness. As another example, 8 identical Lognormals with the lower shape parameter 0.9 and location parameter 5 are randomly added. Order of magnitude for shape value of 0.9 is only about 2.5, and POM is about 1000. The shape of resultant data set in 35,000 simulation runs appears very much Normal-like, easily confirming the prediction of the CLT. The Average of resultant added data is 1778, and which is near the Median value of 1642, indicating that resultant data is almost not skewed at all. Digits are decisively non-Benford, coming at {62.5, 24.0, 4.5, 0.9, 0.4, 0.6, 1.3, 2.4, 3.5}, with an extremely high SSD value of 1353.

Hence, the table in Figure 20 - with its moderate 1.36 average shape parameter for these 8 added Lognormals - should be viewed as a compromise between two extremes, namely between the very high POM (2.2 shape parameter) and the very low POM (0.9 shape parameter). The moderate choice for the shape parameters for Figure 20 was made for pedagogical purposes.

If one still doubts the validity of the general tendency or principle demonstrated here with these three distinct schemes of Lognormal addition processes above, then two more Monte Carlo simulations with k/x distributions would convince even the most avowed skeptic of the validity of the principle. The first scheme is of k/x defined over (10, 100) with its low 10 POM value. Here we randomly add 6 of them, namely: k/x + k/x + k/x + k/x + k/x + k/x. First digits are {29.9, 55.1, 13.9, 0.7, 0.0, 0.0, 0.0, 0.1, 0.2} with an extremely high 1675.2 SSD value. The Average is 234.7 while the Median is 231.8, and the curve of resultant additions is extremely similar to the Normal. Benfordness of the original k/x distribution was completely ruined by the process of these random additions. In sharp contrast, the second scheme of k/x defined over (1, 1000000), with its large 1000000 POM value, is highly resistant to CLT for a long while, not showing any resemblance to the Normal (although it would eventually crumble after many more such random additions). Here as well we randomly add 6 of them, namely: k/x + k/x + k/x + k/x + k/x + k/x. First digits are {28.7, 13.8, 10.7, 9.9, 8.5, 7.9, 7.4, 6.8, 6.5} and SSD value is 30.7, which is fairly low. The Average is 432,976, and significantly larger than the Median value of 293,828, while the curve of resultant additions is not similar to the Normal at all.

This resonates (in some very vague sense) as being in harmony with the new definition of logarithmic-ness given in Kossovsky (2014) chapter 136 page 583, where k/x had to be defined over a huge range in order to qualified as being truly logarithmic with respect to all number system bases and with respect to all bin schemes, being resistant to all such base and bin changes, and thus conserving its logarithmic status. One should not get carried away too far with this analogy, because the two issues are fundamentally different; here we deal with temporary and illusionary resistance to CLT, which would eventually crumble no matter how high order of magnitude happened to be, while the definition of logarithmic-ness is about the conservation of k/x logarithmic status under base and bin changes, and which is another matter.

The moral of the story is that '**not all perfectly logarithmic data sets are created equal**'! Those with very large order of magnitude are by far more resistant to detrimental addition processes in the system, while those with just sufficient order of magnitude for a [nearly] perfect Benford behavior are much more vulnerable to random additions. A data set with say 500,000 points generated via simulations from the Lognormal with 3.5 shape parameter and whatever location, is by far superior and much more resistant to random additions and CLT than a comparable set with say 1.2 shape parameter. Both data sets measure equally for all practical purposes in their current logarithmic status and sizes, yet the former is much more resistant and stable as compared with the latter.

# [12]  Numerical Example of the Achilles' Heel of the Central Limit Theorem

This unique addition game at the casino necessitates having 5 especially-handcrafted 13-sided dice. The faces of these large dice are not marked with usual sequence of the natural integers 1 to 13 as customary with most dice games, but rather each of the 5 dice has the 13 chosen arbitrary numbers of {**3, 7, 12, 19, 24, 35, 42, 76, 92, 176, 331, 564, 978**} written on its sides. The game involves the throwing of these five large dice simultaneously and then the adding of these five values to arrive at the sum upon which a lot of money is bet.

The smallest possible sum here is (3) + (3) + (3) + (3) + (3) = 15.
The biggest possible sum here is (978) + (978) + (978) + (978) + (978) = 4890.

Examples of other possible game occurrences are:

(42) + (42) + (42) +  (7) + (92) =  225
(76) + (978) + (92) + (35) + (24) = 1205
(92) +  (3) + (564) + (24) +  (3) =  686
 (3) + (12) + (19) + (35) +  (7) =   76

The result of 10,000 Monte Carlo computer simulations of this game is depicted in the histogram of Figure 21. The smallest sum in the simulations was 23. The biggest sum in the simulations was 3674. The histogram in Figure 21 depicts almost the entire set of resultant sums, from 20 to 3036, in even steps of 232-length bin width, incorporating 9,964 games, and leaving out only 36 games with extremely big sums of over 3036, not shown in the histogram for lack of space and better visualization.

The histogram shows a marked decline with a tail falling to the right, and nearly consistently skewed pattern favoring the small. Even though only 5 additions are performed here, still, on the face of it, this seems very surprising and quite odd, given the expectation of symmetric and Normal-like histograms for almost all types of addition processes, as predicated by the Central Limit Theorem. Yet, since CLT's Achilles' heel predicts that for skewed variables with high order of magnitude, CLT can manifest itself with difficulties, and only after a truly large number of additions of the random variable, it follows that if our set of 13 numbers on the dice is highly skewed and if it is also of high order of magnitude, then these empirical Monte Carlo result obtained here are certainly expected and reasonable. The crux of the matter is then to examine the quantitative structure of the set of the 13 numbers on the dice, and to determine whether or not it is highly skewed and with high order of magnitude.

Indeed, the set of 13 numbers on the dice is clearly biased against the big and in favor of the small, namely skewed. In order to be able to visualize the skewness of the values of the dice, all 13 values are plotted along the horizontal axis as depicted in Figure 22. The plot actually distorts the scale by stretching it on the left for small values and compressing it on the right for big values, for brevity and better visualization. In spite of giving the big such advantage over the small, the visual message coming out of the plot is that the small is numerous (condensed) and the big is rare (sparse); i.e. that these 13 numbers are skewed. In terms of the non-modern definition of skewness, namely (average − median)/(standard deviation), here the average of the 13 dice numbers is 181.5, and it is much greater than the median which is only 42.0.

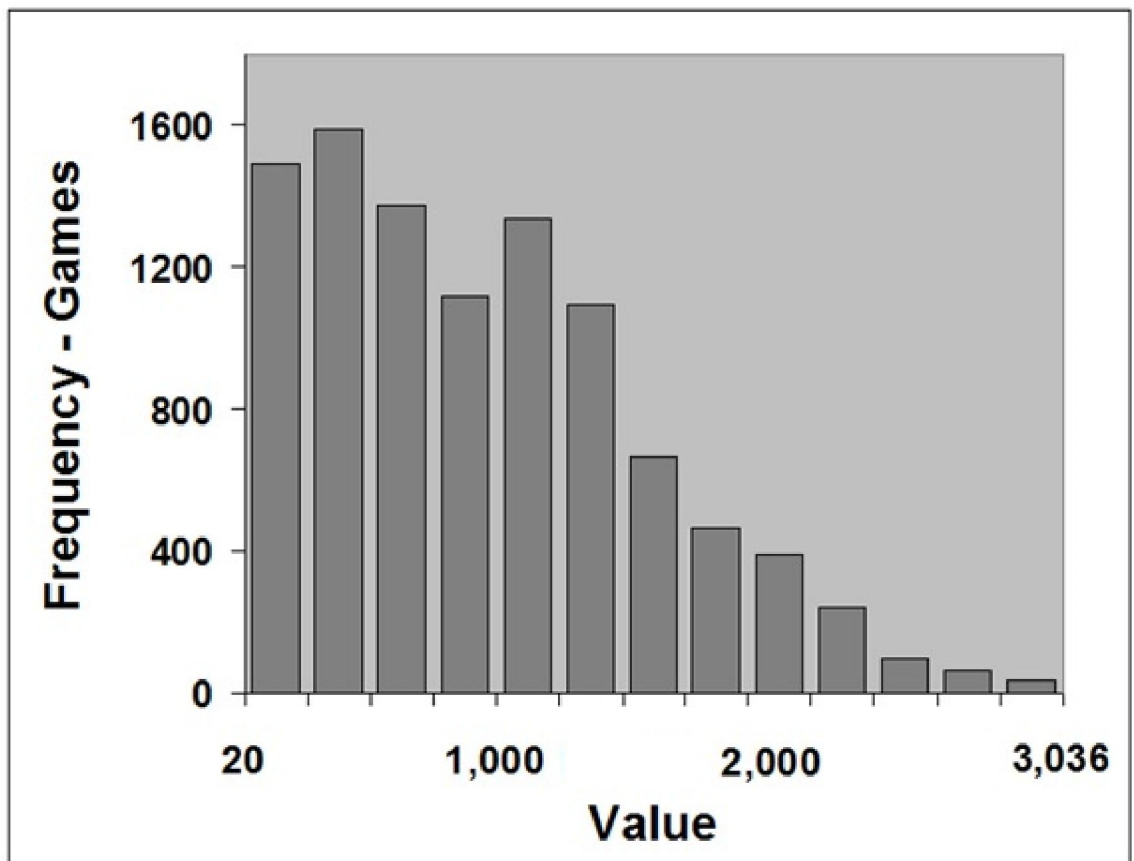

**Figure 21:** Histogram of 10,000 Simulations of the Five 13-Sided Dice Addition Game

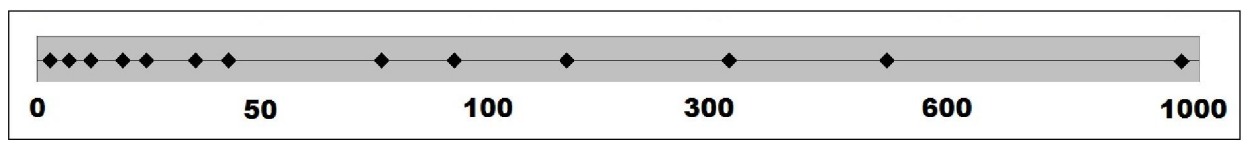

**Figure 22:** Higher Concentration for Small and Diluted Spread for Big for the 13-Sided Dice

In addition, the set of the 13 numbers on the dice is indeed with high order of magnitude. POM is calculated as (978)/(3) = 326, and this value is quite high; especially for example in comparison with a normal 6-sided die of only (6)/(1) = 6 POM value.

Since the set of 13 numbers on the dice is highly skewed and of high order of magnitude, and since only 5 dice are added, the CLT's Achilles' heel clearly explains why CLT could not be manifested here. It would be necessary to build at least about 15 such 13-sided dice for the casino and add 15 dice faces in the game, in order to manifest the CLT well enough. Indeed, computer simulations of 10,000 random additions of 15 such dice (with the same skewed values of high POM above) resulted in nearly symmetrical histogram, resembling the Normal Distribution a great deal, and where the medium came out as the most numerous size by far.

# [13] Tugs of War between Multiplication and Addition

In order to further convince ourselves that these two benign-looking arithmetical operations, namely addition and multiplication, do occasionally go to war with each other over digital and quantitative configurations, another more vivid demonstration of their potential conflicts in the context of Benford's Law is given by the following 6 distinct random processes, where U represents the Uniform(5, 33), and where all Us are identical.

U + U
U*U + U*U
U*U*U + U*U*U
U*U*U*U + U*U*U*U
U*U*U*U*U + U*U*U*U*U
U*U*U*U*U*U + U*U*U*U*U*U

As opposed to adding Lognormals as in Figure 20 where the number of addends is increasing one by one in the main expression, and where the number of multiplicands (hidden within each Lognormal) is fixed, here the number of multiplicands is increasing while the number of addends is fixed at 2. Six Monte Carlo computer Simulations with 35,000 runs each are performed separately at each stage. One should notice carefully that for one realization of say Uniform(5, 33)*Uniform(5, 33), two separate and independent realizations are needed from Uniform(5, 33), followed by the multiplication of these two realizations; and that we are not referring to the simple squaring of a single realization such as in Uniform(5, 33)$^2$.

Simulation results are as follows:

U + U              {  6.3, 19.2, 31.5, 26.8, 13.9,  2.3,  0.0,  0.0, 0.0}  SSD = 1360.0
UU + UU            {21.1,  5.7,  9.4, 11.4, 12.4, 11.8, 10.5,  9.4, 8.2}  SSD =  334.4
UUU + UUU          {42.4, 17.1,  6.7,  4.6,  5.2,  5.6,  6.3,  6.1, 5.9}  SSD =  222.5
UUUU + UUUU        {30.8, 22.4, 14.7,  9.5,  6.7,  5.1,  4.1,  3.5, 3.3}  SSD =   38.9
UUUUU + UUUUU      {26.0, 17.7, 13.9, 11.2,  8.9,  7.2,  6.0,  5.1, 4.1}  SSD =   22.3
UUUUUU + UUUUUU    {30.2, 16.2, 11.5,  9.8,  8.3,  7.1,  6.4,  5.6, 5.0}  SSD =    4.3

Finding addition sleeping at the wheel, being fixed at 2, multiplication then strives hard to win by gradually increasing the number of multiplicands. Finally multiplication is able to achieve the Benford condition, where it manifests itself 6 times within the last expression versus only 2 manifestations of addition. Feeling overconfident now and conceited, multiplication then challenges addition to a new tug of war, willing to face 3 addends – given that it is already in possession of 6 multiplicands:

UUUUUU + UUUUUU + UUUUUU
{35.8, 16.4, 9.6, 7.3, 6.5, 6.6, 6.5, 5.8, 5.5}   SSD = 51.9

Upon seeing the severe setback in digital configuration where SSD is now over fifty, multiplication deeply regrets its previous offer, and it then asks addition for permission to use two more multiplicands, achieving a total of 8 multiplicands versus these 3 addends:

UUUUUUUU + UUUUUUUU + UUUUUUUU
{27.2, 18.5, 13.7, 10.4, 8.5, 6.7, 5.6, 4.8, 4.6}   SSD = 12.1

Multiplication is now quite satisfied with this latest improvement in Benfordness, while addition is not much disturbed upon seeing that digits are getting back fairly close to Benford. Addition knows that ultimately CLT is on its side, and that it would eventually win the war between it and multiplication in due time, and especially when POM is low or when addend distributions are symmetrical and Normality is achieved fairly quickly. Nonetheless, addition wishes to make one last stand in order to flaunt its clout and demonstrate its ability to control the situation, and it expands addends by just one more, namely a total of 4 addends versus these 8 multiplicands:

UUUUUUUU + UUUUUUUU + UUUUUUUU + UUUUUUUU
{25.3, 15.5, 13.7, 11.3, 9.5, 7.9, 6.3, 5.5, 4.8}   SSD = 35.8

Indeed, by adding one more addend to the process, addition was able to move it decisively away from Benford, managing to increase SSD from 12.1 to 35.8.

Where would all these actions and reactions, attacks and counterattacks; that endless tit-for-tat war of attrition between addition and multiplication lead to? Without a doubt, the CLT is decisively on the side of addition, guaranteeing its eventual triumph. Surely multiplication could win some battles in the short run, but it would lose the war in the long run as resultant data finally becomes Normally distributed and digital configuration turning decisively non-Benford, and regardless of the form and defined range of added distributions.

# [14]  The Effects of Random Addition Processes on POM

Let us examine theoretically how random addition processes of an <u>identical</u> random variable affect resultant POM. The random process involves repeatedly adding an identical random variable, over and over again, N times.

**In theory** - POM should not change at all under random additions of N identical random variables X, Y, Z, and so forth. This is so theoretically since the lowest possible value (the worst scenario) for the sum is $MIN_x + MIN_y + MIN_z + \ldots$ (N times), and since X = Y = Z and so forth, $MIN_x = MIN_y = MIN_z$ = (generic MIN), hence the above sum is simply MIN*N. In the same vein, theoretically, the highest possible value (the best scenario) is $MAX_x + MAX_y + MAX_z + \ldots$ (N times), or simply MAX*N, hence $POM_{SUM}$ = (MAX*N)/(MIN*N) = (MAX)/(MIN), namely the same value as the POM of the identical variable being repeatedly added. More succinctly:

## $POM_{SUM} = POM_{THE\ IDENTICAL\ VARIABLE}$

**In reality** - when dealing with numerous additions, this almost never occurs in computer simulations, because it is exceedingly rare that the highest added value is gotten by aiming at all the various maximums simultaneously; and it is exceedingly rare that the lowest added value is gotten by aiming at all the various minimums simultaneously. The maximum added value in simulations is lower than the theoretical; and the minimum added value in simulations is higher than the theoretical; all of which implies that POM of added variables is somewhat lower in actual simulations than the theoretical POM value. Surely for an enormous number of simulation runs, numbering in the billions or trillions say, theoretical POM can be achieved on the computer, and certainly in the abstract, in the limit, as the number of runs goes to infinity. This is why in Figure 20 successive values of POM are monotonically decreasing, as it becomes increasing more difficult for the process to aim at all the maximums and all the minimums simultaneously while more and more additions are applied. Hence in extreme generality for identical variables: the more variables we randomly add, the lower is the resultant POM obtained for the computer simulated addition process - not numbering in the trillions.

As an example, successive random additions of Uniform(3, 17) distributions are run in 8 stages with 35,000 simulation each. In the first stage only one Uniform is simulated 35,000 times without any addition involved. In the 2nd stage two such Uniforms are added randomly 35,000 times, and so forth. The simulations yield:

| # of Vars | 1 | 2 | 3 | 4 | 5 | 6 | 7 | 8 |
|---|---|---|---|---|---|---|---|---|
| Min Simul. | 3.1 | 6.1 | 9.6 | 14.2 | 19.5 | 23.1 | 29.8 | 35.5 |
| Max Simul. | 16.9 | 33.8 | 50.7 | 65.5 | 79.8 | 95.0 | 109.6 | 123.1 |
| POM Simul. | 5.47 | 5.54 | 5.28 | 4.60 | 4.10 | 4.11 | 3.68 | 3.47 |
| POM Theo. | 5.66 | 5.66 | 5.66 | 5.66 | 5.66 | 5.66 | 5.66 | 5.66 |

Theoretical POM here is always 17/3 = 5.66, for any number of random additions whatsoever. Clearly, as can be seen in the above table, actual POM gradually and very slowly decreases from 5.47 just below the theoretical level of 5.66, and finally falling to 3.47 on the 8th additive stage. This is the result of running only 35,000 simulations at each stage, instead of billions or trillions simulations which would yield results very close to the theoretical expectations of POM.
For only two random additions, very little discrepancy is seen between actual and theoretical, namely 5.66 – 5.54 = 0.12, while for say four random additions, much larger discrepancy is seen between actual and theoretical, namely 5.66 – 4.60 = 1.06. This is so since there is a good chance of obtaining in two additions actual POM value of say (16.98 + 16.95)/(3.01 + 3.02), and such fortunate possibility could occur at times, yet, the chance of obtaining in four additions actual POM value of say (16.98 + 16.95 + 16.98 + 16.99)/(3.01 + 3.02 + 3.03 + 3.01) is much smaller, and such fortunate possibility where all 8 extremes are quite close to the top and bottom values of 17 and 3 is quite rare. Put differently, for 2 additions only 4 rare coincidences are needed; while for 4 additions 8 rare coincidences are needed and this scenario is even less likely. The same reasoning applies to <u>multiplication processes</u> with even more dramatic effects, and where ratio of actual POM to theoretical POM decreases even more sharply and rapidly.

As opposed to the theoretical analysis of POM under multiplication processes where the two variables PDF(x) and PDF(y) were allowed to be considered as distinct, here for the analysis of POM under addition processes, without any loss of generality, the scope is restricted initially to identical variables.

Little reflection is needed to realize that this dichotomy [distinct vs. identical] does not restrict the general principles demonstrated here a great deal, and the purpose of considering initially only identical variables for addition processes is merely to ease the analysis and to be able to arrive at one correct conclusion here by being able to cancel the term N in the numerator and in the denominator. Insisting on adding two variables of distinct types and ranges would lead to:

For the random addition process $X + Y$:
The best scenario is adding $MAX_X + MAX_Y$
The worst scenario is adding $MIN_X + MIN_Y$
These two extreme scenarios represent the max and min values for the addition process itself.

$POM_{X+Y} = [ MAX_X + MAX_Y ] / [ MIN_X + MIN_Y ] \neq [2*MAX] / [2*MIN]$
And this inequality is due to the fact that the variables are distinct, spanning distinct ranges.

Yet, in general, the principle that POM of the sum is about as modest as POM of one typical variable still holds true approximately. Allowing variables to be distinct, yet postulating for convenience that the maximum of the minimums is significantly lower than the minimum of the maximums, and that no single variable is much different than the rest; that no single variable overwhelms results by having a range much different than the rest, we can then obtain the approximate (and quite vague) result:

$POM_{AVERAGE} \equiv (1/N)*[ (MAX_X/MIN_X) + (MAX_Y/MIN_Y) + (MAX_Z/MIN_Z) + \ldots ]$

$POM_{SUM} = [ MAX_X + MAX_Y + MAX_Z + \ldots ] / [ MIN_X + MIN_Y + MIN_Z + \ldots ]$

$POM_{SUM} = (1/N)*[ MAX_X + MAX_Y + MAX_Z + \ldots ] / (1/N)*[ MIN_X + MIN_Y + MIN_Z + \ldots ]$

$POM_{SUM} = [ \text{Average Maximum} ] / [ \text{Average Minimum} ] \approx POM_{AVERAGE}$

The main motivation in asserting that $POM_{SUM} \approx POM_{AVERAGE}$, or that at least the two are not much different from each other when no single variable is much different than the rest, is to establish the idea that POM is not being increased due to the act of adding, and that it is still just about at the POM level of any one (average) distribution.

Some concrete numerical examples might assure us that the general line of reasoning above is correct. The following four numerical analyses are of the theoretical highest POM possibilities, or of actual computer simulations with billions or trillions runs. Only one addition process is considered in all of these four examples, namely the random addition process of all the six variables together as in SUM = A + B + C + D + E + F, without stages and without partial additions as was done earlier in the example of Uniform(3, 17).

| Variable | A | B | C | D | E | F | **SUM** |
|---|---|---|---|---|---|---|---|
| Max | 43 | 54 | 13 | 77 | 11 | 15 | **213** |
| Min | 3 | 4 | 3 | 8 | 3 | 4 | **25** |
| POM | 14.3 | 13.5 | 4.3 | 9.6 | 3.7 | 3.8 | **8.5** |

Here POM of the addition process is 8.5, and it fits rather nicely within the six POM values of the individual distinct variables. Value 8.5 is very close to the average POM calculated as 8.2.

| Variable | A | B | C | D | E | F | **SUM** |
|---|---|---|---|---|---|---|---|
| Max | 43 | 54 | 957 | 77 | 11 | 15 | **1157** |
| Min | 3 | 4 | 3 | 8 | 3 | 4 | **25** |
| POM | 14.3 | 13.5 | 319.0 | 9.6 | 3.7 | 3.8 | **46.3** |

Here POM of the addition process is 46.3, which is way over POM of variables A, B, D, E, F, but way below POM of variable C. Here the range for variable C is markedly different than the ranges of the others. Value 46.3 is also quite different from the average POM calculated as 60.6.

| Variable | A | B | C | D | E | F | **SUM** |
|---|---|---|---|---|---|---|---|
| Max | 43 | 54 | 346978 | 77 | 11 | 15 | **347178** |
| Min | 3 | 4 | 5678 | 8 | 3 | 4 | **5700** |
| POM | 14.3 | 13.5 | 61.1 | 9.6 | 3.7 | 3.8 | **60.9** |

Here POM of the addition process is 60.9, which is way over POM of variables A, B, D, E, F, and it closely follows POM of variable C. Value 60.9 is also very different from the average POM calculated as 17.7. Here the huge range for variable C dominates all the others variables.

| Variable | A | B | C | D | E | F | **SUM** |
|---|---|---|---|---|---|---|---|
| Max | 5 | 12 | 54 | 310 | 2144 | 15335 | **17860** |
| Min | 3 | 3 | 3 | 3 | 3 | 3 | **18** |
| POM | 1.7 | 4.0 | 18.0 | 103.3 | 714.7 | 5111.7 | **992.2** |

Here POM of the addition process is 992.2, which is over POM of variables A, B, C, D, E, but well below POM of variable F. Here the ranges steadily increase as focus shifts from A to F. In addition, the 992.2222 POM value for the Sum of the 6 variables is equal almost exactly to the average POM calculated as 992.2333.

These four numerical examples hint at the possibility that POM of the addition process is always less than the largest POM among all the individual variables, namely that it never exceeds the maximum POM in the set of variables, hence the inequality $POM_{SUM} \leq POM_{MAX}$.

Let us prove this in extreme generality, without postulating anything about the variables whatsoever, allowing them to assume any distribution form, and to be defined on all sorts of distinct ranges. [Note: the general relationship $MAX_X = POM_X * MIN_X$ shall be applied.]

$$POM_{SUM} = \frac{MAX_A + MAX_B + MAX_C + MAX_D + \text{etc.}}{MIN_A + MIN_B + MIN_C + MIN_D + \text{etc.}}$$

$$\mathbf{POM_{SUM}} = \frac{POM_A*MIN_A + POM_B*MIN_B + POM_C*MIN_C + POM_D*MIN_D + \text{etc.}}{MIN_A + MIN_B + MIN_C + MIN_D + \text{etc.}}$$

Let us artificially create an expression of $POM_{MAX}$ with a similar format:

$$POM_{MAX} = POM_{MAX} * \frac{MIN_A + MIN_B + MIN_C + MIN_D + \text{etc.}}{MIN_A + MIN_B + MIN_C + MIN_D + \text{etc.}}$$

$$\mathbf{POM_{MAX}} = \frac{POM_{MAX}*MIN_A + POM_{MAX}*MIN_B + POM_{MAX}*MIN_C + POM_{MAX}*MIN_D + \text{etc.}}{MIN_A + MIN_B + MIN_C + MIN_D + \text{etc.}}$$

By definition, $POM_X \leq POM_{MAX}$ for each X variable, hence when comparing the right hand side of the last expression for **POM_SUM** with the right hand side of the last expression for **POM_MAX**, it follows that:

## $POM_{SUM} \leq POM_{MAX}$

For the 10 by 10 multiplication table from our elementary school days converted into an addition table, POM of each {1, 2, 3, ... ,10} addend is 10/1, namely 10, while POM of the entire addition table is 20/2, or 10. This result nicely complies with the expression for identical variables $POM_{SUM} = POM_{THE\ IDENTICAL\ VARIABLE}$. This result also complies with the general expression $POM_{SUM} \leq POM_{MAX}$.

# [15] Summary of the Effects of Arithmetical Processes on Resultant Data

Let us demonstrate in another way the sharp contrast between random multiplication processes and random addition processes. We randomly add four Uniform variables and then contrast all this by randomly multiplying the same set of four Uniform variables. The four variables are: Uniform(6, 75), Uniform(3, 37), Uniform(5, 55), and Uniform(2, 35). We shall give theoretical POM values, as well as actual POM values from computer simulations with 35,000 runs each.

Addition process:     Uniform(6, 75) + Uniform(3, 37) + Uniform(5, 55) + Uniform(2, 35)
Multiplication process:   Uniform(6, 75)*Uniform(3, 37)*Uniform(5, 55)*Uniform(2, 35)

Addition process:          Theoretical POM =  13         Actual POM = 9
Multiplication process:    Theoretical POM = 29677    Actual POM = 6626

Addition process:        1st Digits – {61.8,  0.1,   0.3, 1.1, 2.7, 4.7, 7.6, 9.8, 12.0}  SSD = 1645
Multiplication process: 1st Digits – {30.3, 18.2, 12.7, 9.5, 8.0, 6.5, 5.4, 4.8,   4.4}   SSD = 0.9

The individual POM values are {75/6, 37/3, 55/5, 35/2}, namely {12.5, 12.333, 11.0, 17.5}.

$POM_{SUM} = [MAX_X + MAX_Y + MAX_Z + \ldots ] / [MIN_X + MIN_Y + MIN_Z + \ldots ]$
$POM_{PRODUCT} = [MAX_X * MAX_Y * MAX_Z * \ldots ] / [MIN_X * MIN_Y * MIN_Z * \ldots ]$

Theoretical $POM_{SUM}$      = (75 + 37 + 55 + 35)/(6 + 3 + 5 + 2) = (202)/(16) = 13
Theoretical $POM_{PRODUCT}$ = (75*37*55*35)/(6*3*5*2) = (5341875)/(180) = 29677

Alternatively, applying the expression $POM_{PRODUCT} = \prod POM_J$ regarding POM values of the four Uniform variables as calculated above, we obtain:
Theoretical $POM_{PRODUCT}$ = (12.5)*(12.333)*(11.0)*(17.5) = 29677

For the process in Figure 20 of repeated additions of eight logarithmic Lognormal distributions, POM is steadily decreasing, thereby weakening the intrinsic Benford configuration at each stage. Yet, a worse calamity is awaiting Frank Benford, namely that skewness is also diminishing at each stage of addition, finally attaining that highly symmetrical curve totally contrary to Benford, namely that of the Normal distribution, as predicated by the CLT. Consequently, additive CLT, in one impulsive instant, ruins what MCLT labored on so hard, and for so long, in building up all these individual skewed Lognormals by way of multiplications.

In sharp contrast, Figure 13 examines what happens to POM and SSD values in a random multiplicative process of four Uniform distributions, demonstrating how it all converges to Benford – a sort of a reversal and the opposite of the random additive process of the eight Lognormals examined in Figure 20 where it all diverges away from Benford.

Let us summarize the effects of random arithmetical processes on resultant data:

Randomly multiplying Uniforms yield Lognormal
Randomly multiplying Normals yield Lognormal
Randomly multiplying Lognormals yield Lognormal
Randomly multiplying Benfords yield Benford
**Randomly multiplying non-Benfords yield Benford**

Randomly adding Lognormals yield Normal
Randomly adding Uniforms yield Normal
Randomly adding Normals yield Normal
Randomly adding non-Benfords yield non-Benford
**Randomly adding Benfords yield non-Benford**

<u>Addition Processes:</u>        Less POM  -  More Symmetry -  CLT  -   Normal  - Anti Benford
<u>Multiplication Processes:</u>  More POM  -  More Skewness -  MCLT - Lognormal - Pro Benford

The typical multiplicative form of the equations in physics, chemistry, astronomy, geology, and practically all other scientific disciplines, as well as those of their many applications and results, almost always leads to the manifestation of Benford's Law in the physical world. This is so because resultant order of magnitude is sufficiently high in almost all cases.

# [16]  Multiplication's Ablity to Explain Numerous Logarithmic Processes

Let us apply what was learnt in this article regarding tugs of war between additions and multiplications to actual stochastic and physical models with respect to Benford behavior.

In Kossovsky (2014) chapter 97 titled "The Remarkable Versatility of Benford's Law" discusses the seemingly unrealistic quest to somehow unite all the diverse physical processes, causes, and explanations of Benford's Law into a singular concept, and to show that that fundamental concept is logarithmic. Well, a partial success in at least a limited measure can be found in the generic arithmetical structure of multiplications as examined in this article.

Each logarithmic process has already or should have its own rigorous mathematical proof, and one should at least initially treat them separately, yet, having an overall vista, namely the understanding of the common thread running across all these processes helps in seeing the entire forest instead of just individual and unconnected trees. Let us list eight well-known logarithmic processes that could all come under the protective umbrella of the generic multiplicative process:

1) Random Linear Combinations and accounting data - Kossovsky (2014) chapters 16 and 46.
2) Random Rock Breaking - Kossovsky (2014) chapter 92.
3) Random Multiplications of Random Variables - chapter 8 of this article.
4) Final speed of particles randomly driven to accelerate – chapter 8 of this article.
5) Final position of particles under constant decelerating force – chapter 8 of this article.
6) Exponential Growth Series of the standard deterministic (fixed) factor.
7) Exponential Random Growth Series having a random factor (random log walk) –
   Kossovsky (2014) chapters 22 and 78.
8) Super Exponential Growth Series - Kossovsky (2014) chapter 99.

Surely, other articles will be written about Benford's Law in the future, giving rigorous mathematical proofs that this or that process relating - directly or indirectly, explicitly or implicitly - to multiplication processes is logarithmic. They would all come under the same generic protective umbrella of multiplications. Surely the mathematicians ought to provide distinct rigorous proof for each case separately, yet, given that multiplication is repetitive and that it overwhelms the system (e.g. not being mixed with additions to a great extent as well as having or generating plenty of order of magnitude), the expectation is then to find a decisive logarithmic behavior. It must be stressed though that there are several other processes and data structures that are logarithmic for reasons fundamentally different from multiplications, such as quantitative partitioning, mixture of distributions, data aggregations, chain of distributions, random planet and star formation models as in Kossovsky (2014) chapter 94, and so forth. For this reason, the quest to unite all the diverse physical processes, causes, and explanations in Benford's Law may as well be illusionary.

# [17] Random Linear Combination is not a Generic Explanation of Benford

Typical **Random Linear Combinations** (RLC) models in Kossovsky (2014) chapters 16 and 46 are of the form: [Price List]*Dice, or of the form: [Price List]*Dice1 + [Price List]*Dice2. The term 'Price List' signifies a very short list of prices with only 6 to 9 items typically, such as {$2.25, $4.75, $7.75, $9.50, $10.25, $35.00}, or {$2.25, $4.75, $35.00}, as well as the list {$2.25, $3.25, $4.75, $7.75, $9.50, $10.25, $25.00, $35.00, $37.00}. The term 'Dice' signifies the standard 6-sided die as the discrete set of {1, 2, 3, 4, 5, 6} with uniform and equal probability for each integer. Dice here represents the number of quantities per item bought by the shopper as shown on the face of the randomly thrown die. RLC is a model for the supermarket's revenue data, or equivalently, a model for the shopper's expense data.

**RLC derives its logarithmic tendency exclusively due to the multiplicative nature involved**, namely, the product of the Price List by the Dice, and often this is so in spite of the additive terms involved, therefore logarithmic behavior here is in harmony and consistent with all that was discussed in this article. **Indeed, RLC is <u>not</u> any generic or new explanation of the phenomenon of Benford's Law in and of itself in any sense;** rather RLC is just another manifestation of the ability of the multiplicative process to serve as a generic and authentic explanation of the Benford phenomenon in many real-life cases.

Surely Price List could be structured logarithmically in and of itself, but that's an exogenous issue. And clearly, for a presumed shopper determined to purchase numerous distinct items, say 4 as in [Price List]*Dice1 + [Price List]*Dice2 + [Price List]*Dice3 + [Price List]*Dice4, revenue data may not be Benford if Price List is with low order of magnitude or if it's symmetrical. This is so because the final bill here is structured more as additions than as multiplications, and such state of affairs might let CLT ruin any chance toward Benfordness by pointing to the Normal as the emerging distribution. Benford behavior could still be found here for the 4-item purchase, or even for 5 or 6 items, and so forth, given that Price List is highly skewed <u>and</u> that its order of magnitude is very large, because in such cases the Achilles' heel of the Central Limit Theorem would likely prevent the process from achieving any rapid convergence to the Normal distribution. One should keep in mind that even in cases where the Price List is totally symmetrical, and even though Dice is indeed symmetrical, their product [Price List]*Dice is skewed and non-symmetrical as in all multiplication processes, and therefore if such skewness is strong enough, and if Price List and Dice combine to yield high order of magnitude for their product, then it might invoke the Achilles' heel of the CLT, prevents the Normal from emerging too rapidly, and leads to Benford behavior even for a shopper buying 4, 5, or 6 items.

The understanding gained in this article helps to explain one result in Random Linear Combinations involving continuous variables [as opposed to discrete Price Lists] that hitherto appeared quite peculiar. This result appears in Kossovsky (2014) chapter 46 pages 182 & 183, where **Uniform(0, UB)*Dice1 + Uniform(0, UB)*Dice2 + Uniform(0, UB)*Dice3** strongly deviates from the logarithmic regardless of the particular value assigned to UB, while **Uniform(0, UB)*Dice** was found to be approximately logarithmic for any UB value.

Clearly, in such tugs of war between multiplication and addition as in the 3 Uniforms model, addition should have the upper hand here via CLT, assuming that Uniform(0, UB) and Dice do not combine to yield very high order of magnitude for their product, or that the skewness of Uniform(0, UB)*Dice is not strong enough, and that therefore the rate of convergence of CLT here is rapid; that the Achilles' heel of the Central Limit Theorem is not relevant here, and thus results are decisively non-logarithmic – as indeed was confirmed empirically there.

On the other hand, the single term Uniform(0, UB)*Dice is exclusively multiplicative in nature; totally lacking any additions; the CLT is not involved here at all; and therefore results are nearly logarithmic.

On another note, most of the discrete examples about RLC given in Kossovsky (2014) chapter 46 come out quite close to the logarithmic – in spite of that even and balanced split between addition and multiplication such as in the model [Price List]*Dice1 + [Price List]*Dice2. The (partial) reason for this is that the Price Lists themselves in that chapter are often quite skewed quantitatively and their orders of magnitude are relatively large; rendering the model resistant to CLT by way of the Achilles' heel of the CLT.

## [18]  Distinct Models of Revenue Data Lead to Varied Digital Results

In light of what was discussed in this article, the claim made on page 191 in Kossovsky (2014), that not only [General Store Shopping] but also [Car] is logarithmic under certain conditions, supposedly because both processes generate enough variability in resultant data, cannot be valid, unless CAR's components are exceedingly skewed and come with some extremely large orders of magnitude, so that the process is resistant to CLT by way of the Achilles' heel of the CLT. Surely, such scenario for Car is highly unlikely. In addition, Car has to overcome the big hurdle of having so many additive terms within its expression, rendering it highly vulnerable to CLT's tenacious drive towards Normality.

[General Store Shopping] = $N_1$*Item1 + $N_2$*Item2 + $N_3$*Item3

[Car] = 4* Brake + 1*Engine (with its many smaller components) + 2*Bumper +

4*Door + 1*Fuel Tank + 3*Mirror + 1*Transmission System + 4*Bearing + 4*wheel +

1*Air Conditioning + 1*Steering System + 4*Tire + 2*Air Bag + etc.

The model of General Store Shopping expressed as the addition of <u>three</u> addends appears to resemble a bit more addition processes than multiplication processes, hence given that order of magnitude of the price list or the catalog is not high enough, or that it is nearly symmetrical, results might be significantly influenced by the CLT and thus Normal-like, resulting in deviation from Benford.

On the other hand, the distribution of General Store Shopping expressed as the addition of only <u>two</u> addends [*as is usually the case in Kossovsky (2014) chapter 46*] resembles multiplication processes and addition processes in equal measure, and thus results might be close to Benford, especially when Price List is highly skewed and its order of magnitude is high.

In summary: Random Linear Combinations and accounting data are normally full of tensions and dramatic tugs of wars between additions and multiplications. Therefore each particular manifestation or model of Random Linear Combinations and accounting data should be judged and evaluated individually according to the relative strength of the competing addition and multiplication forces and their abilities to exert influence upon the system, as well as according to the quantitative configuration of the price list involved.

## [19]  Revenue Data and the Achilles' Heel of the CLT

The typical bill for a shopper at a large supermarket or at a big retail store purchasing several different items is a mixture of random combinations of multiplications and additions. A single bill for one typical shopper in a large supermarket or at a big retail store may read as follows:

**3***($**2.75** *bread*) + **5***($**2.50** *tuna*) + **2***($**7.99** *cheese*) = ($**36.73** *total bill*)

Since the typical shopper buys only very <u>few distinct items</u>, it follows that very few additions are involved, and therefore CLT cannot manifest itself well, especially if the price list is skewed and of high order of magnitude. Empirical examinations of typical price lists and catalogs of several well-known and big retail shops, supermarkets, and large companies show that they are indeed highly skewed in favor of the small, and that their orders of magnitude are quite large. Therefore, the typical bill should be highly skewed in favor of the small and Benford, just as the underlying price lists and catalogs are, and this conclusion is indeed in perfect agreement with empirical examinations of real-life bills in revenue data sets.

Clearly, for a particularly presumed shopper determined to purchase <u>numerous distinct items</u>, revenue data may not be skewed in favor of the small and may not be Benford. But even for these shoppers, the Achilles' heel of the CLT could still strongly retard convergence to the Normal distribution if the price list is highly skewed and of high order of magnitude.

For more modest retail stores and smaller supermarkets having price lists and catalogs without much skewness and with relatively modest order of magnitude, revenue data may deviate somewhat from Benford. Decisive deviation from Benford in revenue data can be found in very small business operations, such as a very small coffee shop at the corner of an insignificant and short street, serving only 3 brands of overpriced and very bitter coffee, as well as just 2 types of tasteless cakes, and having a very small clientele.

# [20] Expense Data Empirical Compliance with Benford's Law

Almost all accounting, financial, and economics-related data – and with only very few and very rare exceptions - are structured in such a way that the small is more numerous than the big, namely having positive skewness, and with the Benford digit configuration.

Empirical examination of actual revenue data at the raw level, detailed bill by bill (as opposed to summaries or ratios) is nearly impossible to obtain due to confidentiality and secrecy issues.

Only Quarterly or Annual Financial Statements are public information that can be readily examined, but Benford and skewness are encountered mostly in the raw detailed values (at the individual bill, receipt, and invoice level), and not really in summary values, aggregations, or ratios of values, such as the numbers typically found in the Financial Statements.

Fortunately for statisticians, some governmental expense data is available at the raw and original level, detailing each and every transaction.

The State Of Oklahoma in the USA provides detailed information at the transaction level for all its vendor payments for the fiscal year 2011. The website for this database is at: https://data.ok.gov/dataset/state-oklahoma-vendor-payments-fiscal-year-2011

These payments reflect disbursements from a state fund for the purchase of goods received, services performed, reimbursements, and payments to other governments.

Although this data set is purely of expenses and costs, yet, it also reflects strongly on revenue data in general. This is so since every entry here for a given expense is surely also one revenue item for some provider, company, or agent, billing the state and charging it for the product or service rendered.

Examination of this data set reveals that the distribution of expense amounts is highly skewed in favor of the small. Figure 23 depicts the histogram of the vast majority of expenses from $0 to $1 million, containing 986,962 items; leaving out only a small minority of 530 very expensive items of over $1 million. The horizontal x-axis scale is mostly a logarithmic one, except for that jump or gap from $0 to $100. Clearly, except for a brief and very gentle rise on the left of the histogram between $0 and $1000, the cheap (considered as small) outnumbers the expensive (considered as big). The temporary and minor reversal of the histogram in the beginning for very low values is quite typical in revenue and expense accounting data, yet, the overall description of relative quantities is decisively in favor of the small.

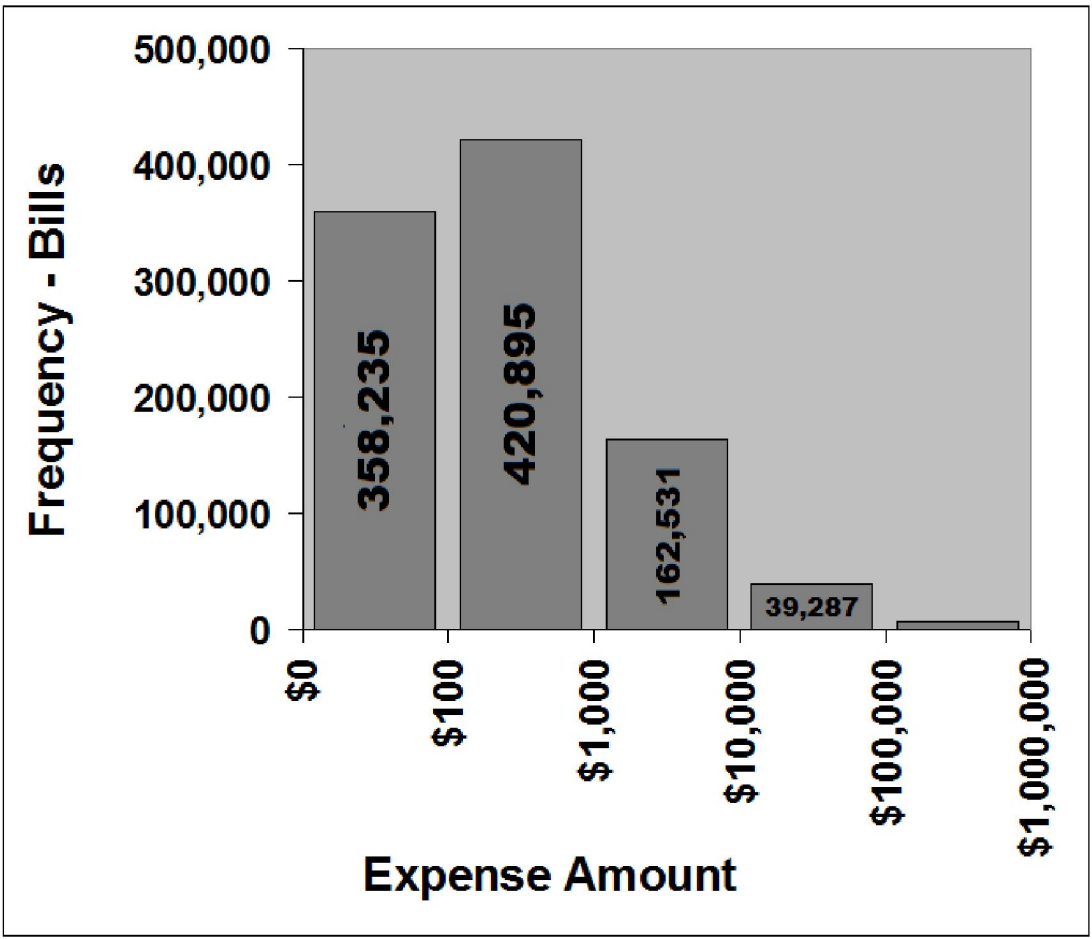

**Figure 23**: Histogram of the State of Oklahoma Expenses is Skewed in Favor of Small

Examining digital configuration for this data set (of the expenses from $0 to $1 million, containing 986,962 items) yields the following results:

First digits for Oklahoma: {29.7, 17.7, 12.1, 9.7, 8.6, 6.6, 6.1, 4.9, 4.5}
Benford's Law 1st Digits : {30.1, 17.6, 12.5, 9.7, 7.9, 6.7, 5.8, 5.1, 4.6}

The extremely low 0.9 SSD value here is a strong indication that the Oklahoma expenses data set complies with Benford's Law almost perfectly. Figure 24 visually demonstrates the strong compliance with Benford's Law for this extremely large data set containing 987,492 bills.

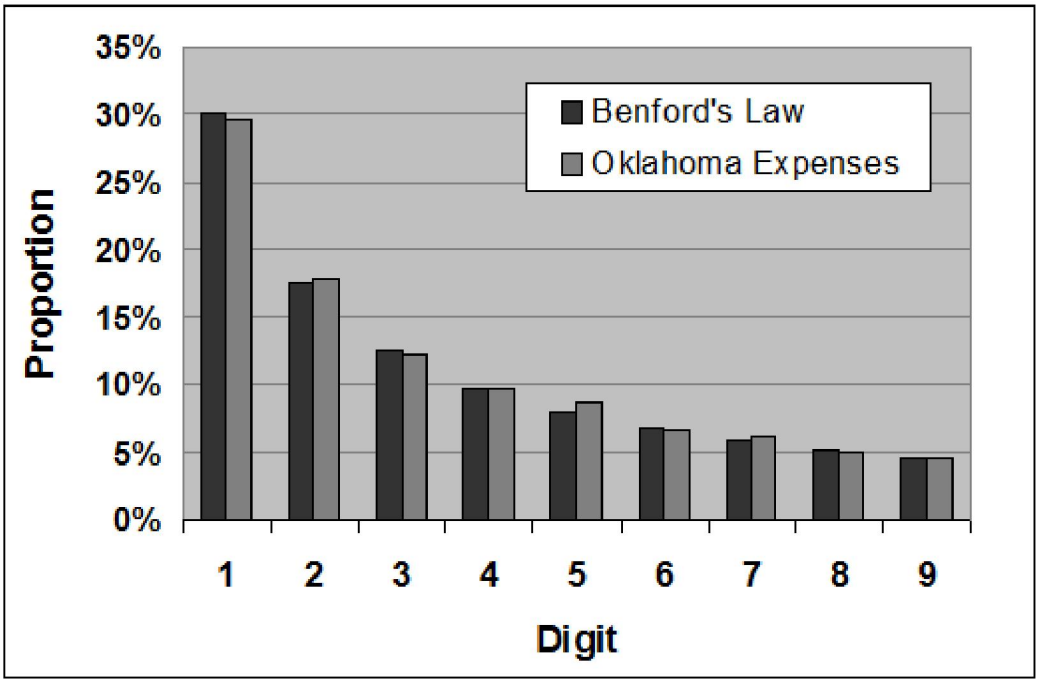

**Figure 24**: Compliance of State of Oklahoma Expenses with Benford's Law

### [21] Molar Mass of Chemical Compounds and the Achilles' Heel of the CLT

Lactose is chosen as a good representative of chemical compounds, demonstrating how molecules are typically formed from elementary atoms in the natural world. The chemical formula for Lactose is **$C_{12}H_{22}O_{11}$** and its molar mass is 342.29648 g/mol.

**12*(12.0107** *Carbon*) + **22*(1.00794** *Hydrogen*) + **11*(15.9994** *Oxygen*) = **342.29648** g/mol

The approach here differs only slightly from the case of revenue data pertaining to a large collection of purchasing bills. The consideration is of the molar mass, namely the weight of the entire molecule, which is simply the sum of the weights of all the small atoms constituting the large molecule.

Empirical examination of the Periodic Table shows that its histogram is nearly symmetrical and that it's only marginally and very slightly skewed in favor of the small. In addition, order of magnitude of the Periodic Table is found to be fairly large.

The main part of the Periodic Table is quantitatively analyzed, considering the most relevant elements, from Hydrogen with atomic number 1, all the way to Radon with atomic number 86 and atomic weight 222.0; and ignoring the less relevant heavier elements with atomic numbers over 86. Figure 25 depicts the histogram of the Periodic Table by the weight variable, from 0 up to 216 (just below Radon) in steps of 27 unified atomic mass units (u) for the bin width.

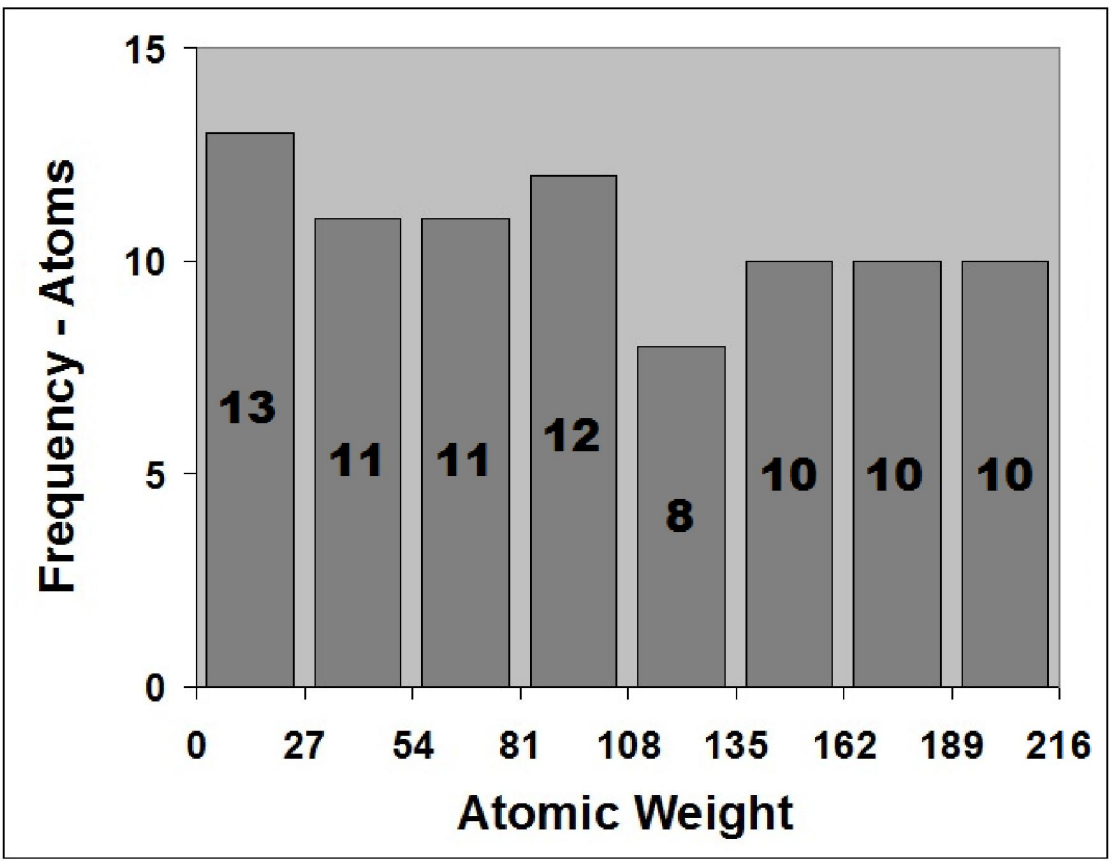

**Figure 25:** Histogram of Atomic Weights in the Periodic Table up to Element Radon

On the face of it, the Central Limit Theorem should guarantee that the weights of molecules are symmetrically distributed; closely mimicking the shape of the Normal Distribution, and that Benford's Law is not obeyed at all. This assumption is based on the observation that the molar mass is structured primarily as an addition process with 3 or 4 typical addends, and even with 5 addends in some cases; while the multiplicative terms in the molar mass expression are with only 2 multiplicands. Yet, this is not the case in this tug of war, and in reality the Achilles' heel of the CLT strongly retards convergence to the Normal distribution. This is so because the multiplicative terms within the expressions are somewhat skewed and of sufficiently high order of magnitude. Indeed, since the Periodic Table itself is of high order of magnitude, these multiplicative terms are of high order of magnitude as well. In addition, these multiplicative terms are also somewhat skewed, even though the Periodic Table is far from being skewed, and this is so since the very act of multiplying - such as in (11)*(Oxygen Mass) – leads to skewness as in all multiplication processes with sufficient order of magnitude (and this act of multiplying also leads to even higher order of magnitude as an extra bonus).

In conclusion: the molar mass of molecules should be skewed in favor of the small and approximately Benford.

# [22]  Molar Mass Empirical Compliance with Benford's Law

A large list of 2175 commonly used and naturally occurring chemical compounds is provided on the website: http://www.convertunits.com/compounds/. Here we consider a large collection of molecules which are relevant to human use, industrial production, human environment, and general availability. The website criterion for the selection of these 2175 molecules does not follow any particular formal procedure, and instead it simply pulls together information from a variety of chemical and scientific sources, as well as utilizing the informal suggestions of chemists and scientists regarding which molecules should be considered as relevant and important for compilation.

The real-life data on molar mass in the above-mentioned large list of 2175 chemical compounds is indeed quantitatively skewed in favor of the small, as predicted theoretically. Figure 26 depicts the histogram up to 900 gram/mole, where small molecules outnumber big molecules in general, except for a brief rise on the very left of the histogram between 0 and 300 where relatively bigger molecules are slightly more numerous than smaller molecules. Such temporary and minor reversal of the histogram in the beginning for very low values is quite typical in many other physical, scientific, financial, and accounting data sets.

Examining compliance with Benford's Law for this set of 2175 chemical compounds we get:

First digits of Molar Mass: {31.9,  25.2,  16.1,  8.4,  5.7,  4.3,  2.9,  3.2,  2.3}
Benford's Law 1st Digits : {30.1,  17.6,  12.5,  9.7,  7.9,  6.7,  5.8,  5.1,  4.6}

Admittedly, the 102.9 SSD value for the Molar Mass data is somewhat high, indicating weak or partial compliance with Benford's Law. Figure 27 depicts the partial compliance of this list of molar mass values with Benford's Law.

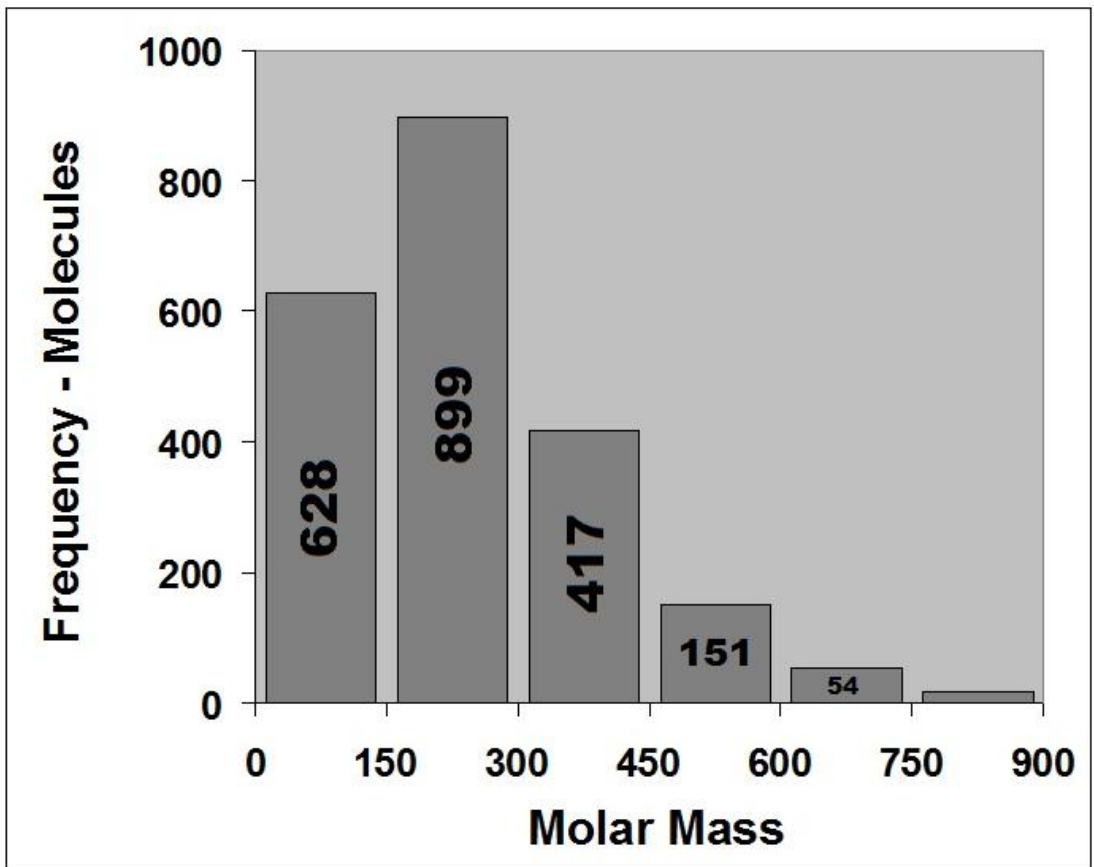

**Figure 26**: Small Molecules Generally Outnumber Big Molecules in the Chemical World

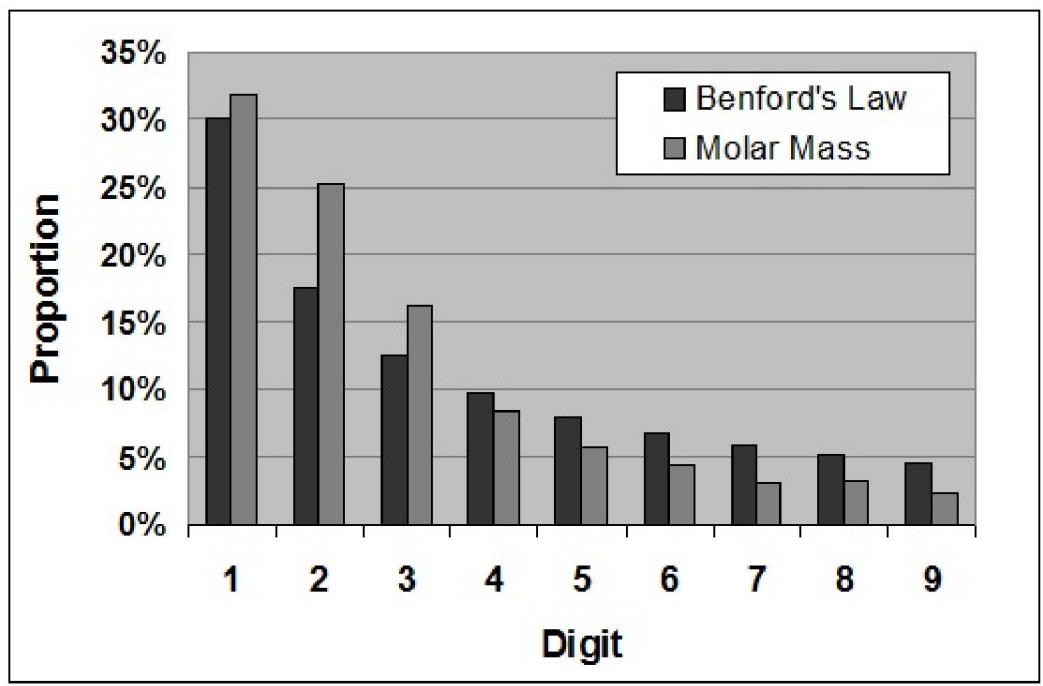

**Figure 27**: Partial Compliance of Molar Mass with Benford's Law

# [23] Applications to Random Consolidations and Fragmentations Models

The model of Random Consolidations and Fragmentations cycles (C&F) represents a classic case where multiplication overcomes addition, thereby winning that life-and-death struggle between these two arithmetical adversaries.

Link to an article by this author on the topic titled "Random Consolidations and Fragmentations Cycles Lead to Benford" can be found at: https://arxiv.org/abs/1505.05235

V is the initial value of each ball (quantity), L is the number of balls in the system, and C is the number of full cycles executed. Schematically the process is described as follows:

1) Initial Set = {V, V, V, …   L times   …V, V, V}
2) Repeat C times:
    (i) Choose one ball at random and split it as in {Uniform(0, 1), 1 – Uniform(0, 1)}.
    (ii) Choose two balls randomly and merge them.
3) Final Set is Benford assuming that many balls L are involved, and that many cycles C are executed, namely that approximately L > 300 or preferably that L > 400, as well as that C > 2*L or preferably that C > 3*L approximately.

In that article, as an example, one particular random trajectory of events is calculated, beginning with five quantitatively identical balls {1, 1, 1, 1, 1} and leading the following evolved five terms after six cycles, where $U_j$ represents a realized value from the Uniform(0, 1):

(U1*U2+2)*U6
(U1*U2+2)*(1-U6)
U1*U3*U5*(1-U2) + (2-U1)*(1-U4)
U1*U3*(1-U2)*(1-U5) + 1 + (2-U1)*U4
U1*(1-U2)*(1-U3)

For this particular random trajectory of events above, out of five expressions, two are sums (minority), and three are products (majority). In fact, since each full cycle yields **3** newly created balls, **2** of which are multiplicative, namely X(U) and X*(1 – U), and **1** of which is an additive, therefore the statistical tendency of the system after numerous cycles is to have approximately 2/3 multiplicative expressions and 1/3 additive expressions.

On the face of it, the existence of additive expressions in about one third of all expressions, does not bode very well for Benford digital configuration; the glass is two-third full and one-third empty, and one doesn't know whether he or she should be happy or should be sad.

Each fragmentation process contributes to the system two products - each with a minimum of 2 multiplicands and possibly more; similarly each consolidation process contributes to the system an additive expression with a minimum of 2 addends and possibly more; and therefore there exists a tug of war here between additions and multiplications with respect to Benford behavior. Remarkably, even though Benford frequently loses numerous C&F battles, yet he wins the war in the long run. One existing feature here that is partially saving the system from deviation from the logarithmic is that on average only about one-third of the expressions are additive; and that even within those expressions there are plenty of arithmetical multiplicative elements involved. Surely there are some additive expressions with 3 addends which might be quite detrimental to Benford, but they are far and few between; and there are even more menacing expressions with 4 addends, but luckily these are even rarer.

The general [theoretical] understanding gained in this article regarding multiplication and addition processes enables us to thoroughly explain the [empirical] strong logarithmic behavior in all Consolidations and Fragmentations models; namely the reason addition effects do not manage to significantly retard multiplication effects. In a nutshell, the C&F process is Benford because it uses the high OOM variable of Uniform(0, 1) which contributes to high order of magnitude, skewness, and thus Benfordness. As a consequence, the C&F process encounters the Achilles' heel of the Central Limit Theorem and additions are not very effective. Order of magnitude of the Uniform(0, 1) calculated as LOG(1/0) is infinite. A more conservative measure of OOM with the exclusion of the edges or 'outliers', calculated as $Q_{90\%}$ divided by $Q_{10\%}$ would yield 0.9/0.1 or simply 9, which is finite. Surely, the C&F model cannot use any low OOM variable such as say Uniform(5, 7), because it needs to break a whole quantity into two fractions, and this can only be achieved via the high OOM variable Uniform(0, 1). Such high OOM values, coupled with the fact that the terms within the additive expressions almost always involve also some multiplications (which are always skewed), guarantee that the Central Limit Theorem is very slow to act here and that its retarded rate of convergence does not manage to even begin to ruin the general multiplicative tendencies of the system. Skewness for these multiplicative terms hiding within the additive expressions is guaranteed by virtue of simply being multiplication. All multiplication processes yield skewed set of values. Since the vast majority of the additive terms here are with only 2, 3, or 4 addends, the Central Limit Theorem does not even begin to manifests itself.

## [24] Application of the Distributive Rule in Random Arithmetical Processes

One should be cautious in interpreting the application of the Distributive Rule in algebra to expressions describing random statistical processes. For example, given the probability distributions A, B, C, D, X; the process of adding A, B, C, D, followed by the multiplication of that sum by X is symbolically written as **X*(A + B + C + D)**, namely that the process is principally a product of two multiplicands, and on the face of it may be quite close to Benford assuming large enough orders of magnitude for the individual distributions. Yet, the application of the Distributive Rule here would yield **X*A + X*B + X*C + X*D**, namely that the process is principally an addition of four addends, and on the face of it may be non-Benford due to CLT, especially if distributions are symmetrical or if they are with low orders of magnitude.

The above interpretation is mistaken. These two expressions above signify an identical statistical process which is multiplicative in nature, in spite of their distinct appearances! It is important to keep in mind that in the latter expression of X*A + X*B + X*C + X*D, all four X symbols denote an identical and singular value of the variable, in other words, one particular realization of X in the stimulations. A truly distinct random process could be gotten by having four distinct and independent realizations of X, such as in the expression $X_1*A + X_2*B + X_3*C + X_4*D$, where it truly resembles more additions than multiplications; most likely be non-Benford; and closely resembling the Normal Distribution. Surely, the Distributive Rule cannot be applied here since in general $X_I \neq X_J$.

## [25] An Historical Note on Simon Newcomb's Ratio of Uniforms Assertion

In light of the new results and better understanding obtained in this article, Simon Newcomb's two-page short article in 1881 now appears to contain not only the correct digital proportion in real life typical data sets (i.e. Benford's Law), but also a remarkable insight into why it should be so in data relating to the physical sciences - namely 'almost' one correct explanation for the existence of the law!

Newcomb writes: *"As natural numbers occur in nature, they are to be considered as the ratios of quantities. Therefore, instead of selecting a number at random, we must select two numbers, and inquire what is the probability that the first significant digit of their ratio is the digit n."*

While the ratio [or product] of just two independent random variables may not be quite exactly "Benford" or "Newcomb", and digital configuration still depends somewhat on the type of distributions involved and especially on their orders of magnitude, yet it could be very close to the logarithmic, as was seen in this article. Clearly, ratios and products belong to the same class of random arithmetical processes, while additions and subtractions belong to another class.

Alex Ely Kossovsky, Jan 1, 2019
akossovsky@gmail.com